\documentclass[12pt,twoside]{amsart}
\usepackage[latin1]{inputenc}
\usepackage{amsmath, amsthm, amscd, amsfonts, amssymb, graphicx}
\usepackage[bookmarksnumbered, plainpages]{hyperref}

\newtheorem{thm}{Theorem}[section]

\newtheorem{lem}[thm]{Lemma}

\newtheorem{defn}[thm]{Definition}

\numberwithin{equation}{section}

\begin{document}

\title{\bf The general Kastler-Kalau-Walze type theorem for the $J$-twist $D_{J}$ of the Dirac operator }
\author{Siyao Liu \hskip 0.4 true cm  Yong Wang$^{*}$}

\thanks{{\scriptsize
\hskip -0.4 true cm \textit{2010 Mathematics Subject Classification:}
53C40; 53C42.
\newline \textit{$^{*}$Corresponding author}}}

\maketitle

\begin{abstract}
 \indent In \cite{LSW1} and \cite{LSW2}, we proved the Kastler-Kalau-Walze type theorem for the $J$-twist $D_{J}$ of the Dirac operator on $3$-dimensional, $4$-dimensional and $6$-dimensional almost product Riemannian spin manifold with boundary.
 In this paper, we generalize our previous conclusions and establish the proof of the general Kastler-Kalau-Walze type theorem for the $J$-twist $D_{J}$ of the Dirac operator on even-dimensional almost product Riemannian spin manifold with boundary.
\end{abstract}

\vskip 0.2 true cm


\pagestyle{myheadings}
\markboth{\rightline {\scriptsize Liu}}
         {\leftline{\scriptsize  A Kastler-Kalau-Walze type theorem for the $J$-twist $D_{J}$ of the Dirac operator}}

\bigskip
\bigskip


\section{ Introduction }
The noncommutative residue was found in \cite{Gu,Wo}. Since noncommutative residues are of great importance to the study of noncommutative geometry, more and more attention has been attached to the study of noncommutative residues.
Connes put forward that the noncommutative residue of the square of the inverse of the Dirac operator was proportioned to the Einstein-Hilbert action, which is called the Kastler-Kalau-Walze theorem now \cite{Co1,Co2}. Kastler gave a brute-force proof of this theorem \cite{Ka}. 
In the same time, Kalau-Walze and Ackermann proved this theorem by using the normal coordinates system and the heat kernel expansion, respectively \cite{KW,Ac}. 
The result of Connes was extended to the higher dimensional case \cite{U}. 
Wang generalized the Connes' results to the case of manifolds with boundary in \cite{Wa1,Wa2}, and proved the Kastler-Kalau-Walze type theorem for the Dirac operator and the signature operator on lower-dimensional manifolds with boundary. 
In \cite{Wa3}, for the Dirac operator, Wang computed $\widetilde{{\rm Wres}}[\pi^+D^{-1}\circ\pi^+D^{-1}]$, in these cases the boundary term vanished.
Most of the operators studied in the literature have the leading symbol $\sqrt{-1}c(\xi)$ \cite{Wa4,Wa5,WW,WWW,Wa}. However, Wu and Wang studied operators with the leading symbol $-\widehat{c}(V)c(\xi)$. In \cite{WW2}, Wu and Wang gave the proof of Kastler-Kalau-Walze type theorems  of the operators $\sqrt{-1}\widehat{c}(V)(d+\delta)$ and $-\sqrt{-1}\widehat{c}(V)(d+\delta)$ on 3,4-dimensional oriented compact manifolds with or without boundary.
On the other hand, some preliminaries and lemmas about the Dirac operator $D$ and the $J$-twist are given in \cite{K}. In \cite{Chen1,Chen2}, the author got estimates on the higher eigenvalues of the Dirac operator on locally reducible Riemannian manifolds by the $J$-twist of the Dirac operator. It can be obtained by simple calculations that the leading symbol of the $J$-twist $D_{J}$ of the Dirac operator is not $\sqrt{-1}c(\xi)$.
Therefore, Liu and Wang proved the Kastler-Kalau-Walze type theorem for the $J$-twist $D_{J}$ of the Dirac operator on $3$-dimensional, $4$-dimensional and $6$-dimensional almost product Riemannian spin manifold with boundary \cite{LSW1,LSW2}.

In \cite{Wa6} and \cite{Wa7}, the authors established the general Kastler-Kalau-Walze type theorems for any dimensional manifolds with boundary which generalizes the results in \cite{Wa4,Wa5}.
{\bf The motivation of this paper} is to prove the general Kastler-Kalau-Walze type theorem for the $J$-twist $D_{J}$ of the Dirac operator on even-dimensional almost product Riemannian spin manifold with boundary.
The novelty of our paper is to the computations related to the negative higher power of $D_{J}$.
Motivated by \cite{LSW1} and \cite{LSW2}, we compute the generalized noncommutative residue $\widetilde{{\rm Wres}}[\pi^+(D_{J}^{-1})\circ\pi^+(D_{J}^{-n+3})]$ on even-dimensional manifold with boundary.
Our main theorem is as follows.

\begin{thm}\label{thm1.1}
Let $M$ be an $n$-dimensional almost product Riemannian spin manifold with the boundary $\partial M$, then we get the following equality:
\begin{align}\label{1.1}
&\widetilde{{\rm Wres}}[\pi^+(D_{J}^{-1})\circ\pi^+(D_{J}^{-n+3})]=\frac{(n-2)\pi^{\frac{n}{2}}}{(\frac{n}{2}-1)!}\int_{M}2^{\frac{n}{2}}\Big(\frac{1}{4}\sum_{i,j=1}^{n}R(J(e_{i}), J(e_{j}), e_{j}, e_{i})\\
&-\frac{1}{2}\sum_{\nu,j=1}^{n}g^{M}(\nabla_{e_{j}}^{L}(J)e_{\nu}, (\nabla^{L}_{e_{\nu}}J)e_{j})-\frac{1}{2}\sum_{\nu,j=1}^{n}g^{M}(J(e_{\nu}), (\nabla^{L}_{e_{j}}(\nabla^{L}_{e_{\nu}}(J)))e_{j}-(\nabla^{L}_{\nabla^{L}_{e_{j}}e_{\nu}}(J))e_{j})\nonumber\\
&-\frac{1}{4}\sum_{\alpha,\nu,j=1}^{n}g^{M}(J(e_{\alpha}), (\nabla^{L}_{e_{\nu}}J)e_{j})g^{M}((\nabla^{L}_{e_{\alpha}}J)e_{j}, J(e_{\nu}))-\frac{1}{4}\sum_{\alpha,\nu,j=1}^{n}g^{M}(J(e_{\alpha}), (\nabla^{L}_{e_{\alpha}}J)e_{j})\nonumber\\
&\cdot g^{M}(J(e_{\nu}), (\nabla^{L}_{e_{\nu}}J)e_{j})+\frac{1}{4}\sum_{\nu,j=1}^{n}g^{M}((\nabla^{L}_{e_{\nu}}J)e_{j}, (\nabla^{L}_{e_{\nu}}J)e_{j}))-\frac{5}{12}s\Big)d{\rm Vol_{M} }+\int_{\partial M}\Big[2^{\frac{4-n}{2}}K\nonumber\\
&\cdot\Big(\frac{\pi(n-2)^{2}}{n^{2}+n-2}\cdot\frac{(n-3)!}{(\frac{n}{2})!(\frac{n}{2}-2)!}+\sum_{i=1}^{n}\langle J(e_{i}), e_{n}\rangle^{2}\cdot\frac{\pi n(n^{2}-8n+12)}{8(n-1)}\cdot\frac{(n-3)!}{(\frac{n}{2}+1)!(\frac{n}{2})!}\nonumber\\
&-\langle J(e_{n}), e_{n}\rangle^{2}\cdot\frac{\pi(n^{2}-3n-2)}{n^{2}+n-2}\cdot\frac{(n-3)!}{(\frac{n}{2})!(\frac{n}{2}-2)!}\Big)-2^{-\frac{2+n}{2}}\sum_{i=1}^{n}g^{M}(J(e_{i}), (\nabla^{L}_{e_{i}}J)e_{n})\nonumber\\
&\cdot\pi(n^{2}-8n+12)\cdot\frac{(n-2)!}{(\frac{n}{2}+1)!(\frac{n}{2}-1)!}\Big]Vol(S^{n-2})d{\rm Vol_{\partial M}},\nonumber
\end{align}
where $n$ is the even number, $s$ is the scalar curvature, $Vol(S^{n-2})$ is the canonical volume of $S^{n-2}$ and $K$ are defined in section $3$.
\end{thm}

A brief description of the organization of this paper is as follows.
In Section 2, this paper will first introduce the basic notions of Boutet de Monvel's calculus and the definition of the noncommutative residue for manifolds with boundary. 
In the next section, we recall some basic facts and formulas about the $J$-twist $D_{J}$ of the Dirac operator and prove the proof of the general Kastler-Kalau-Walze type theorem for the $J$-twist $D_{J}$ of the Dirac operator on even-dimensional almost product Riemannian spin manifold with boundary.


\vskip 1 true cm

\section{ Boutet de Monvel's calculus and the definition of the noncommutative residue }

In this section, we explain the basic notions of  Boutet de Monvel's calculus and the definition of the noncommutative residue for manifolds with boundary that will be used throughout the paper. For the details, see Ref.\cite{Wa3}.

Let $M$ be an n-dimensional compact oriented manifold with the boundary $\partial M$.
We assume that the metric $g^{TM}$ on $M$ has the following form near the boundary,
\begin{equation}
g^{M}=\frac{1}{h(x_{n})}g^{\partial M}+dx _{n}^{2},
\end{equation}
where $g^{\partial M}$ is the metric on $\partial M$ and $h(x_n)\in C^{\infty}([0, 1)):=\{\widehat{h}|_{[0,1)}|\widehat{h}\in C^{\infty}((-\varepsilon,1))\}$ for some $\varepsilon>0$ and $h(x_n)$ satisfies $h(x_n)>0$, $h(0)=1$, where $x_n$ denotes the normal directional coordinate. 
Let $U\subset M$ be a collar neighborhood of $\partial M$ which is diffeomorphic with $\partial M\times [0,1)$. 
By the definition of $h(x_n)\in C^{\infty}([0,1))$ and $h(x_n)>0$, there exists $\widehat{h}\in C^{\infty}((-\varepsilon,1))$ such that $\widehat{h}|_{[0,1)}=h$ and $\widehat{h}>0$ for some sufficiently small $\varepsilon>0$. 
Then there exists a metric $g'$ on $\widetilde{M}=M\bigcup_{\partial M}\partial M\times(-\varepsilon,0]$ which has the form on $U\bigcup_{\partial M}\partial M\times (-\varepsilon,0 ]$
\begin{equation}
\label{b2}
g'=\frac{1}{\widehat{h}(x_{n})}g^{\partial M}+dx _{n}^{2} ,
\end{equation}
such that $g'|_{M}=g$. 
We fix a metric $g'$ on the $\widetilde{M}$ such that $g'|_{M}=g$.

Let the Fourier transformation $F'$ be
\begin{equation*}
F':L^2({\bf R}_t)\rightarrow L^2({\bf R}_v);~F'(u)(v)=\int_\mathbb{R} e^{-ivt}u(t)dt
\end{equation*}
and let
\begin{equation*}
r^{+}:C^\infty ({\bf R})\rightarrow C^\infty (\widetilde{{\bf R}^+});~ f\rightarrow f|\widetilde{{\bf R}^+};~
\widetilde{{\bf R}^+}=\{x\geq0;x\in {\bf R}\}.
\end{equation*}
We define $H^+=F'(\Phi(\widetilde{{\bf R}^+}));~ H^-_0=F'(\Phi(\widetilde{{\bf R}^-}))$ which satisfies $H^+\bot H^-_0$, where $\Phi(\widetilde{{\bf R}^+}) =r^+\Phi({\bf R})$, $\Phi(\widetilde{{\bf R}^-}) =r^-\Phi({\bf R})$ and $\Phi({\bf R})$ denotes the Schwartz space. We have the following property: $h\in H^+~$ (resp. $H^-_0$) if and only if $h\in C^\infty({\bf R})$ which has an analytic extension to the lower (resp. upper) complex half-plane $\{{\rm Im}\xi<0\}$ (resp. $\{{\rm Im}\xi>0\})$ such that for all nonnegative integer $l$,
\begin{equation*}
\frac{d^{l}h}{d\xi^l}(\xi)\sim\sum^{\infty}_{k=1}\frac{d^l}{d\xi^l}(\frac{c_k}{\xi^k}),
\end{equation*}
as $|\xi|\rightarrow +\infty,{\rm Im}\xi\leq0$ (resp. ${\rm Im}\xi\geq0)$ and where $c_k\in{\bf C}$ are some constants.
 
Let $H'$ be the space of all polynomials and $H^-=H^-_0\bigoplus H';~H=H^+\bigoplus H^-.$ 
Denote by $\pi^+$ (resp. $\pi^-$) the projection on $H^+$ (resp. $H^-$). Let $\widetilde H=\{$rational functions having no poles on the real axis$\}$. 
Then on $\widetilde{H}$,
\begin{equation}
\label{b3}
\pi^+h(\xi_0)=\frac{1}{2\pi i}\lim_{u\rightarrow 0^{-}}\int_{\Gamma^+}\frac{h(\xi)}{\xi_0+iu-\xi}d\xi,
\end{equation}
where $\Gamma^+$ is a Jordan closed curve included ${\rm Im}(\xi)>0$ surrounding all the singularities of $h$ in the upper half-plane and
$\xi_0\in {\bf R}$. 
In our computations, we only compute $\pi^+h$ for $h$ in $\widetilde{H}$. Similarly, define $\pi'$ on $\widetilde{H}$,
\begin{equation}
\label{b4}
\pi'h=\frac{1}{2\pi}\int_{\Gamma^+}h(\xi)d\xi.
\end{equation}
So $\pi'(H^-)=0$. 
For $h\in H\bigcap L^1({\bf R})$, $\pi'h=\frac{1}{2\pi}\int_{{\bf R}}h(v)dv$ and for $h\in H^+\bigcap L^1({\bf R})$, $\pi'h=0$.

An operator of order $m\in {\bf Z}$ and type $d$ is a matrix\\
$$\widetilde{A}=\left(\begin{array}{lcr}
  \pi^+P+G  & K  \\
   T  &  \widetilde{S}
\end{array}\right):
\begin{array}{cc}
\   C^{\infty}(M,E_1)\\
 \   \bigoplus\\
 \   C^{\infty}(\partial{M},F_1)
\end{array}
\longrightarrow
\begin{array}{cc}
\   C^{\infty}(M,E_2)\\
\   \bigoplus\\
 \   C^{\infty}(\partial{M},F_2)
\end{array},
$$
where $M$ is a manifold with boundary $\partial M$ and $E_1,E_2$~ (resp. $F_1,F_2$) are vector bundles over $M~$ (resp. $\partial M$).
Here, $P:C^{\infty}_0(\Omega,\overline {E_1})\rightarrow C^{\infty}(\Omega,\overline {E_2})$ is a classical pseudo-differential operator of order $m$ on $\Omega$, where $\Omega$ is a collar neighborhood of $M$ and $\overline{E_i}|M=E_i~(i=1,2)$. 
$P$ has an extension: ${\mathcal{E'}}(\Omega,\overline {E_1})\rightarrow{\mathcal{D'}}(\Omega,\overline {E_2})$, where ${\mathcal{E'}}(\Omega,\overline {E_1})~({\mathcal{D'}}(\Omega,\overline{E_2}))$ is the dual space of $C^{\infty}(\Omega,\overline{E_1})~(C^{\infty}_0(\Omega,\overline {E_2}))$. 
Let $e^+:C^{\infty}(M,{E_1})\rightarrow{\mathcal{E'}}(\Omega,\overline{E_1})$ denote extension by zero from $M$ to $\Omega$ and $r^+:{\mathcal{D'}}(\Omega,\overline{E_2})\rightarrow{\mathcal{D'}}(\Omega, {E_2})$ denote the restriction from $\Omega$ to $X$, then define
$$\pi^+P=r^+Pe^+:C^{\infty}(M,{E_1})\rightarrow {\mathcal{D'}}(\Omega,{E_2}).$$

In addition, $P$ is supposed to have the transmission property; this means that, for all $j,k,\alpha$, the homogeneous component $p_j$ of order $j$ in the asymptotic
expansion of the symbol $p$ of $P$ in local coordinates near the boundary satisfies:\\
$$\partial^k_{x_n}\partial^\alpha_{\xi'}p_j(x',0,0,+1)=
(-1)^{j-|\alpha|}\partial^k_{x_n}\partial^\alpha_{\xi'}p_j(x',0,0,-1),$$
then $\pi^+P:C^{\infty}(M,{E_1})\rightarrow C^{\infty}(M,{E_2})$. Let $G$, $T$ be respectively the singular Green operator and the trace operator of order $m$ and type $d$. 
Let $K$ be a potential operator and $S$ be a classical pseudo-differential operator of order $m$ along the boundary. Denote by $B^{m,d}$ the collection of all operators of
order $m$ and type $d$,  and $\mathcal{B}$ is the union over all $m$ and $d$.

Recall that $B^{m,d}$ is a Fr\'{e}chet space. 
The composition of the above operator matrices yields a continuous map:
$B^{m,d}\times B^{m',d'}\rightarrow B^{m+m',{\rm max}\{m'+d,d'\}}.$ 
Write $$\widetilde{A}=\left(\begin{array}{lcr}
 \pi^+P+G  & K \\
 T  &  \widetilde{S}
\end{array}\right)
\in B^{m,d},
 \widetilde{A}'=\left(\begin{array}{lcr}
\pi^+P'+G'  & K'  \\
 T'  &  \widetilde{S}'
\end{array} \right)
\in B^{m',d'}.$$
The composition $\widetilde{A}\widetilde{A}'$ is obtained by multiplication of the matrices (For more details see \cite{Sc}). 
For example $\pi^+P\circ G'$ and $G\circ G'$ are singular Green
operators of type $d'$ and
$$\pi^+P\circ\pi^+P'=\pi^+(PP')+L(P,P').$$
Here $PP'$ is the usual composition of pseudo-differential operators and $L(P,P')$ called leftover term is a singular Green operator of type $m'+d$. 
For our case, $P,P'$ are classical pseudo-differential operators, in other words $\pi^+P\in \mathcal{B}^{\infty}$ and $\pi^+P'\in \mathcal{B}^{\infty}$ .

Let $M$ be an $n$-dimensional compact oriented manifold with the boundary $\partial M$.
Denote by $\mathcal{B}$ the Boutet de Monvel's algebra. We recall that the main theorem in \cite{FGLS}.
\begin{thm}\label{th:32}{\rm\cite{FGLS}}{\bf(Fedosov-Golse-Leichtnam-Schrohe)}
Let $M$ and $\partial M$ be connected, ${\rm dim}M=n\geq3$, and let $\widetilde{S}$ (resp. $\widetilde{S}'$) be the unit sphere about $\xi$ (resp. $\xi'$) and $\sigma(\xi)$ (resp. $\sigma(\xi')$) be the corresponding canonical $n-1$ (resp. $(n-2)$) volume form.
Set $\widetilde{A}=\left(\begin{array}{lcr}\pi^+P+G &   K \\
T &  \widetilde{S}    \end{array}\right)$ $\in \mathcal{B}$ , and denote by $p$, $b$ and $s$ the local symbols of $P, G$ and $\widetilde{S}$ respectively.
Define:
\begin{align}
{\rm{\widetilde{Wres}}}(\widetilde{A})=&\int_X\int_{\bf \widetilde{ S}}{\rm{tr}}_E\left[p_{-n}(x,\xi)\right]\sigma(\xi)dx \nonumber\\
&+2\pi\int_ {\partial X}\int_{\bf \widetilde{S}'}\left\{{\rm tr}_E\left[({\rm{tr}}b_{-n})(x',\xi')\right]+{\rm{tr}}
_F\left[s_{1-n}(x',\xi')\right]\right\}\sigma(\xi')dx',
\end{align}
where ${\rm{\widetilde{Wres}}}$ denotes the noncommutative residue of an operator in the Boutet de Monvel's algebra.
Then~~ a) ${\rm \widetilde{Wres}}([\widetilde{A},B])=0 $, for any
$\widetilde{A},B\in\mathcal{B}$;~~ b) It is the unique continuous trace on
$\mathcal{B}/\mathcal{B}^{-\infty}$.
\end{thm}

Let $p_1$, $p_2$ be non-negative integers and $p_1+p_2\leq n$.
\begin{defn}{\rm\cite{Wa3}}
Lower dimensional volumes of spin manifolds with boundary are defined by
 \begin{equation}
{\rm Vol}^{(p_1,p_2)}_nM:= \widetilde{{\rm Wres}}[\pi^+D^{-p_1}\circ\pi^+D^{-p_2}].
\end{equation}
\end{defn}
Denote by $\sigma_{l}(D)$ the $l$-order symbol of the Dirac operator $D$.
By \cite{Wa3}, we get
\begin{align}\label{b1}
\widetilde{{\rm Wres}}[\pi^+D^{-p_1}\circ\pi^+D^{-p_2}]=\int_M\int_{|\xi|=1}{\rm
trace}_{S(TM)}[\sigma_{-n}(D^{-p_1-p_2})]\sigma(\xi)dx+\int_{\partial M}\Phi
\end{align}
and
\begin{align}\label{b2}
\Phi&=\int_{|\xi'|=1}\int^{+\infty}_{-\infty}\sum^{\infty}_{j, k=0}\sum\frac{(-i)^{|\alpha|+j+k+1}}{\alpha!(j+k+1)!}
\times {\rm trace}_{S(TM)}[\partial^j_{x_n}\partial^\alpha_{\xi'}\partial^k_{\xi_n}\sigma^+_{r}(D^{-p_1})(x',0,\xi',\xi_n)
\\
&\times\partial^\alpha_{x'}\partial^{j+1}_{\xi_n}\partial^k_{x_n}\sigma_{l}(D^{-p_2})(x',0,\xi',\xi_n)]d\xi_n\sigma(\xi')dx',\nonumber
\end{align}
where the sum is taken over $r+l-k-|\alpha|-j-1=-n,~~r\leq -p_1,l\leq -p_2$.

Since $[\sigma_{-n}(D^{-p_1-p_2})]|_M$ has the same expression as $\sigma_{-n}(D^{-p_1-p_2})$ in the case of manifolds without boundary, so locally we can compute the first term by \cite{Ka,KW,Wa3,Po}.

\section{ The noncommutative residue $\widetilde{{\rm Wres}}[\pi^+(D_{J}^{-1})\circ\pi^+(D_{J}^{-n+3})]$ on even-dimensional manifolds with boundary }

We give some definitions and basic notions which we will use in this paper.
Let $M$ be an $n$-dimensional ($n\geq 6$) oriented compact Riemannian manifold with a Riemannian metric $g^{M}$.
We recall that the Dirac operator $D$ is locally given as follows:
\begin{equation}
D=\sum_{i, j=1}^{n}g^{ij}c(\partial_{i})\nabla_{\partial_{j}}^{S}=\sum_{i=1}^{n}c(e_{i})\nabla_{e_{i}}^{S},
\end{equation}
where $\partial_{i}$  is a natural local frame on $TM,$ $(g^{ij})_{1\leq i,j\leq n}$ is the inverse matrix associated to the metric matrix  $(g_{ij})_{1\leq i,j\leq n}$ on $M$ and $c(e_{i})$ be the Clifford action which satisfies the relation
\begin{align}
&c(e_{i})c(e_{j})+c(e_{j})c(e_{i})=-2g^{M}(e_{i}, e_{j})=-2\delta_i^j,
\end{align}
\begin{align}
&\nabla_{\partial_{j}}^{S}=\partial_{i}+\sigma_{i}
\end{align}
and
\begin{align}
&\sigma_{i}=\frac{1}{4}\sum_{j, k=1}^{n}\langle \nabla_{\partial_{i}}^{L}e_{j}, e_{k}\rangle c(e_{j})c(e_{k}).
\end{align}

Let $J$ be a $(1, 1)$-tensor field on $(M, g^M)$ such that $J^2=\texttt{id},$
\begin{align}
&g^M(J(X), J(Y))=g^M(X, Y),
\end{align}
for all vector fields $X,Y\in \Gamma(TM).$ Here $\texttt{id}$ stands for the identity map. $(M, g^M, J)$ is an almost product Riemannian manifold. We can also define on almost product Riemannian spin manifold the following $J$-twist $D_{J}$ of the Dirac operator $D$ by
\begin{align}
&D_{J}:=\sum_{i=1}^{n}c(e_{i})\nabla^{S}_{J(e_{i})}=\sum_{i=1}^{n}c[J(e_{i})]\nabla^{S}_{e_{i}}.
\end{align}

Let $\xi=\sum_{k}\xi_{j}dx_{j},$  $\nabla^L_{\partial_{i}}\partial_{j}=\sum_{k}\Gamma_{ij}^{k}\partial_{k},$  we denote that
\begin{align}
 \Gamma^{k}=g^{ij}\Gamma_{ij}^{k};\ \sigma^{j}=g^{ij}\sigma_{i};\ \partial^{j}=g^{ij}\partial_{i},\nonumber
\end{align}
where $\Gamma_{ij}^{k}$ is the Christoffel coefficient of Levi-Civita connection $\nabla^{L}$.

In this section, we want to compute the noncommutative residue $\widetilde{{\rm Wres}}[\pi^+(D_{J}^{-1})\circ\pi^+(D_{J}^{-n+3})]$ on $n$-dimensional almost product Riemannian spin manifold with the boundary.
By (\ref{b1}) and (\ref{b2}), we can rewrite the noncommutative residue as
\begin{align}\label{c1}
\widetilde{{\rm Wres}}[\pi^+({D}_{J}^{-1})\circ\pi^+({D}_{J}^{-n+3})]=\int_M\int_{|\xi|=1}{\rm
trace}_{S(TM)}[\sigma_{-n}({D}_{J}^{-n+2})]\sigma(\xi)dx+\int_{\partial M}\Phi
\end{align}
and
\begin{align}\label{c2}
\Phi&=\int_{|\xi'|=1}\int^{+\infty}_{-\infty}\sum^{\infty}_{j, k=0}\sum\frac{(-i)^{|\alpha|+j+k+1}}{\alpha!(j+k+1)!}
\times {\rm trace}_{S(TM)}[\partial^j_{x_n}\partial^\alpha_{\xi'}\partial^k_{\xi_n}\sigma^+_{r}({D}_{J}^{-1})(x',0,\xi',\xi_n)
\\
&\times\partial^\alpha_{x'}\partial^{j+1}_{\xi_n}\partial^k_{x_n}\sigma_{l}({D}_{J}^{-n+3})(x',0,\xi',\xi_n)]d\xi_n\sigma(\xi')dx',\nonumber
\end{align}
where the sum is taken over $r+l-k-|\alpha|-j-1=-n,~~r\leq -1,l\leq -n+3$.

By \cite{LSW1} and \cite{LSW2}, we have
\begin{thm}\cite{LSW1,LSW2}
If $M$ is an $n$-dimensional almost product Riemannian spin manifold without boundary, we have the following:
\begin{align}
{\rm Wres}(D_{J}^{-n+2})
=\frac{(n-2)\pi^{\frac{n}{2}}}{(\frac{n}{2}-1)!}\int_{M}2^{\frac{n}{2}}\Big(&\frac{1}{4}\sum_{i,j=1}^{n}R(J(e_{i}), J(e_{j}), e_{j}, e_{i})
-\frac{1}{2}\sum_{\nu,j=1}^{n}g^{M}(\nabla_{e_{j}}^{L}(J)e_{\nu}, (\nabla^{L}_{e_{\nu}}J)e_{j})\\
&-\frac{1}{2}\sum_{\nu,j=1}^{n}g^{M}(J(e_{\nu}), (\nabla^{L}_{e_{j}}(\nabla^{L}_{e_{\nu}}(J)))e_{j}-(\nabla^{L}_{\nabla^{L}_{e_{j}}e_{\nu}}(J))e_{j})\nonumber\\
&-\frac{1}{4}\sum_{\alpha,\nu,j=1}^{n}g^{M}(J(e_{\alpha}), (\nabla^{L}_{e_{\nu}}J)e_{j})g^{M}((\nabla^{L}_{e_{\alpha}}J)e_{j}, J(e_{\nu}))\nonumber\\
&-\frac{1}{4}\sum_{\alpha,\nu,j=1}^{n}g^{M}(J(e_{\alpha}), (\nabla^{L}_{e_{\alpha}}J)e_{j})g^{M}(J(e_{\nu}), (\nabla^{L}_{e_{\nu}}J)e_{j})\nonumber\\
&+\frac{1}{4}\sum_{\nu,j=1}^{n}g^{M}((\nabla^{L}_{e_{\nu}}J)e_{j}, (\nabla^{L}_{e_{\nu}}J)e_{j}))-\frac{5}{12}s\Big)d{\rm Vol_{M} }.\nonumber
\end{align}
\end{thm}

Locally we can use Theorem 3.1 to compute the interior term of (\ref{c1}), then
\begin{align}
&\int_M\int_{|\xi|=1}{\rm
trace}_{S(TM)}[\sigma_{-n}({D}_{J}^{-n+2})]\sigma(\xi)dx\\
&=\frac{(n-2)\pi^{\frac{n}{2}}}{(\frac{n}{2}-1)!}\int_{M}2^{\frac{n}{2}}\Big(\frac{1}{4}\sum_{i,j=1}^{n}R(J(e_{i}), J(e_{j}), e_{j}, e_{i})
-\frac{1}{2}\sum_{\nu,j=1}^{n}g^{M}(\nabla_{e_{j}}^{L}(J)e_{\nu}, (\nabla^{L}_{e_{\nu}}J)e_{j})\nonumber\\
&-\frac{1}{2}\sum_{\nu,j=1}^{n}g^{M}(J(e_{\nu}), (\nabla^{L}_{e_{j}}(\nabla^{L}_{e_{\nu}}(J)))e_{j}-(\nabla^{L}_{\nabla^{L}_{e_{j}}e_{\nu}}(J))e_{j})\nonumber\\
&-\frac{1}{4}\sum_{\alpha,\nu,j=1}^{n}g^{M}(J(e_{\alpha}), (\nabla^{L}_{e_{\nu}}J)e_{j})g^{M}((\nabla^{L}_{e_{\alpha}}J)e_{j}, J(e_{\nu}))\nonumber
\end{align}
\begin{align}
&-\frac{1}{4}\sum_{\alpha,\nu,j=1}^{n}g^{M}(J(e_{\alpha}), (\nabla^{L}_{e_{\alpha}}J)e_{j})g^{M}(J(e_{\nu}), (\nabla^{L}_{e_{\nu}}J)e_{j})\nonumber\\
&+\frac{1}{4}\sum_{\nu,j=1}^{n}g^{M}((\nabla^{L}_{e_{\nu}}J)e_{j}, (\nabla^{L}_{e_{\nu}}J)e_{j}))-\frac{5}{12}s\Big)d{\rm Vol_{M} }.\nonumber
\end{align}

In the following, we will compute $\int_{\partial M}\Phi$.
As shown in \cite{LSW1} and \cite{LSW2}, we see that
\begin{lem}\cite{LSW1}
The following identities hold:
\begin{align}
\sigma_1({D}_{J})&=ic[J(\xi)];\\
\sigma_0({D}_{J})&=
-\frac{1}{4}\sum_{i,j,k=1}^{n}\omega_{j,k}(e_i)c[J(e_i)]c(e_j)c(e_k),
\end{align}
where $(\omega_{j,k})$ denotes the connection matrix of Levi-Civita connection $\nabla^{L}$.
\end{lem}
\begin{lem}\cite{LSW1}
The following identities hold:
\begin{align}
\sigma_{-1}({D}_{J}^{-1})&=\frac{ic[J(\xi)]}{|\xi|^2};\\
\sigma_{-2}({D}_{J}^{-1})&=\frac{c[J(\xi)]\sigma_{0}({D}_{J})c[J(\xi)]}{|\xi|^4}+\frac{c[J(\xi)]}{|\xi|^6}\sum_ {j=1}^{n} c[J(dx_j)]
\Big[\partial_{x_j}(c[J(\xi)])|\xi|^2-c[J(\xi)]\partial_{x_j}(|\xi|^2)\Big].
\end{align}
\end{lem}
\begin{lem}\cite{LSW2} 
The following identities hold:
\begin{align}
\sigma_3({D}_{J}^{3})&=ic[J(\xi)]|\xi|^{2};\\
\sigma_2({D}_{J}^{3})&=\sum^{n}_{i,j,l=1}c[J(dx_{l})]\partial_{l}(g^{ij})\xi_{i}\xi_{j}+c[J(\xi)](2\sigma^k-\Gamma^k)\xi_{k}\\
&-\sum^{n}_{\alpha=1}c[J(\xi)]c[J(e_{\alpha})]c[(\nabla^{L}_{e_{\alpha}}J)(\xi^{*})]-\frac{1}{4}|\xi|^2\sum^{n}_{s,t,l=1}\omega_{s,t}(e_{l})c[J(e_{l})]c(e_{s})c(e_{t}),\nonumber
\end{align}
where $\xi^{*}=\sum^{n}_{\beta=1}\langle e_{\beta}, \xi\rangle e_{\beta}.$
\end{lem}

According to (2.13) in \cite{LSW1}, we have
\begin{lem} 
The following identities hold:
\begin{align}
\sigma_2({D}_{J}^{2})&=|\xi|^{2};\\
\sigma_1({D}_{J}^{2})&=i(-2\sigma^k+\Gamma^k)\xi_{k}+i\sum^{n}_{\alpha=1}c[J(e_{\alpha})]c[(\nabla^{L}_{e_{\alpha}}J)(\xi^{*})].
\end{align}
\end{lem}

Write
 \begin{eqnarray}
D_x^{\alpha}&=(-i)^{|\alpha|}\partial_x^{\alpha};
~\sigma({D}_{J}^{2})=p_{2}+p_{1}+p_{0};
~\sigma({D}_{J}^{-2})=\sum^{\infty}_{j=2}q_{-j}.
\end{eqnarray}
By the composition formula of pseudodifferential operators, we have
\begin{align}
1=\sigma({D}_{J}^2\circ {D}_{J}^{-2})&=
\sum_{\alpha}\frac{1}{\alpha!}\partial^{\alpha}_{\xi}
[\sigma({D}_{J}^2)]{{D}}^{\alpha}_{x}
[\sigma({D}_{J}^{-2})] \\
&=(p_2+p_1+p_0)(q_{-2}+q_{-3}+\cdots)\nonumber\\
&+\sum_j(\partial_{\xi_{j}}p_2+\partial_{\xi_j}p_1+\partial_{\xi_{j}}p_0)
(D_{x_j}q_{-2}+D_{x_j}q_{-3}+\cdots) \nonumber\\
&=p_2q_{-2}+(p_2q_{-3}+p_1q_{-2}+\sum_j\partial_{\xi_j}p_2D_{x_j}q_{-2})+\cdots,\nonumber
\end{align}
so
\begin{equation}
q_{-2}=p_2^{-1};~q_{-3}=-p_2^{-1}[p_1p_2^{-1}+\sum_j\partial_{\xi_j}p_2D_{x_j}(p_2^{-1})].
\end{equation}
\begin{lem} The following identities hold:
\begin{align}
\sigma_{-2}({D}_{J}^{-2})&=|\xi|^{-2};\\
\sigma_{-3}({D}_{J}^{-2})&=-i|\xi|^{-4}(-2\sigma^k+\Gamma^k)\xi_{k}-i|\xi|^{-4}\sum^{n}_{\alpha=1}c[J(e_{\alpha})]c[(\nabla^{L}_{e_{\alpha}}J)(\xi^{*})]\\
&-i|\xi|^{-6}\sum^{n}_{j=1}2\xi^j\xi_\alpha\xi_\beta\partial_{x_{j}}(g^{\alpha\beta}).\nonumber
\end{align}
\end{lem}

Since $\sigma_{-2}({D}_{J}^{-2})=|\xi|^{-2}$, it follows that $\sigma_{-n}({D}_{J}^{-n})=|\xi|^{-n}$ and 
\begin{align}
\sigma_{-n+3}({D}_{J}^{-n+3})&=\sigma_{-n+3}({D}_{J}^{-n}\cdot{D}_{J}^{3})\\
&=|\xi|^{-n}\cdot ic[J(\xi)]|\xi|^{2}=ic[J(\xi)]|\xi|^{-n+2}.\nonumber
\end{align}

From \cite{Wa7}, we similarly can obtain
\begin{align}
\sigma_{-n+2}({D}_{J}^{-n+3})&=\sigma_{-n+2}({D}_{J}^{-n+2}\cdot{D}_{J})\\
&=\Bigg\{\sum_{|\alpha|=0}^{+\infty}\frac{1}{\alpha!}\partial_{\xi}^{\alpha}[\sigma({D}_{J}^{-n+2})]D_{x}^{\alpha}[\sigma({D}_{J})]\Bigg\}_{-n+2}\nonumber\\
&=\sigma_{-n+2}({D}_{J}^{-n+2})\sigma_{0}({D}_{J})-i\sum_{j=1}^{n}\partial_{\xi_{j}}[\sigma_{-n+2}({D}_{J}^{-n+2})]\partial_{x_{j}}[\sigma_{1}({D}_{J})]\nonumber\\
&+\Bigg\{\frac{n-2}{2}\sigma^{(-\frac{n}{2}+2)}_{2}({D}_{J}^{2})\sigma_{-3}({D}_{J}^{-2})-i\sum_{k=0}^{\frac{n}{2}-3}\sum_{\mu=1}^{n}\partial_{\xi_{\mu}}\sigma^{(-\frac{n}{2}+k+2)}_{2}({D}_{J}^{2})\nonumber
\end{align}
\begin{align}
&\times \partial_{x_{\mu}}\sigma^{-1}_{2}({D}_{J}^{2})(\sigma_{2}({D}_{J}^{2}))^{(-k)}\Bigg\}\sigma_{1}({D}_{J}).\nonumber
\end{align}

Since $\Phi$ is a global form on $\partial M,$ so for any fixed point $x_0\in\partial M$, we can choose the normal coordinates $U$ of $x_0$ in $\partial M$ (not in $M$) and compute $\Phi(x_0)$ in the coordinates $\widetilde{U}=U\times [0,1)\subset M$ and the metric $\frac{1}{h(x_n)}g^{\partial M}+dx_n^2.$ The dual metric of $g^M$ on $\widetilde{U}$ is ${h(x_n)}g^{\partial M}+dx_n^2.$ 
Write $g^M_{ij}=g^M(\frac{\partial}{\partial x_i},\frac{\partial}{\partial x_j});~ g_M^{ij}=g^M(dx_i,dx_j)$, then

\begin{equation}
[g^M_{ij}]= \left[\begin{array}{lcr}
  \frac{1}{h(x_n)}[g_{ij}^{\partial M}]  & 0  \\
   0  &  1
\end{array}\right];~~~
[g_M^{ij}]= \left[\begin{array}{lcr}
  h(x_n)[g^{ij}_{\partial M}]  & 0  \\
   0  &  1
\end{array}\right]
\end{equation}
and
\begin{equation}
\partial_{x_s}g_{ij}^{\partial M}(x_0)=0, 1\leq i,j\leq n-1; ~~~g_{ij}^M(x_0)=\delta_{i}^{j}.
\end{equation}

$\{e_1, \cdots, e_n\}$ be an orthonormal frame field in $U$ about $g^{\partial M}$ which is parallel along geodesics and $e_i(x_0)=\frac{\partial}{\partial{x_i}}(x_0).$ We review the following three lemmas.
\begin{lem}{\rm \cite{Wa3}}\label{le:32}
With the metric $g^{M}$ on $M$ near the boundary
\begin{eqnarray}
\partial_{x_j}(|\xi|_{g^M}^2)(x_0)&=&\left\{
       \begin{array}{c}
        0,  ~~~~~~~~~~ ~~~~~~~~~~ ~~~~~~~~~~~~~{\rm if }~j<n, \\[2pt]
       h'(0)|\xi'|^{2}_{g^{\partial M}},~~~~~~~~~~~~~~~~~~~~{\rm if }~j=n;
       \end{array}
    \right. \\
\partial_{x_j}[c(\xi)](x_0)&=&\left\{
       \begin{array}{c}
      0,  ~~~~~~~~~~ ~~~~~~~~~~ ~~~~~~~~~~~~~{\rm if }~j<n,\\[2pt]
\partial x_{n}(c(\xi'))(x_{0}), ~~~~~~~~~~~~~~~~~{\rm if }~j=n,
       \end{array}
    \right.
\end{eqnarray}
where $\xi=\xi'+\xi_{n}dx_{n}$.
\end{lem}
\begin{lem}{\rm \cite{Wa3}}\label{le:32}With the metric $g^{M}$ on $M$ near the boundary
\begin{align}
\omega_{s,t}(e_i)(x_0)&=\left\{
       \begin{array}{c}
        \omega_{n,i}(e_i)(x_0)=\frac{1}{2}h'(0),  ~~~~~~~~~~ ~~~~~~~~~~~{\rm if }~s=n,t=i,i<n; \\[2pt]
       \omega_{i,n}(e_i)(x_0)=-\frac{1}{2}h'(0),~~~~~~~~~~~~~~~~~~~{\rm if }~s=i,t=n,i<n;\\[2pt]
    \omega_{s,t}(e_i)(x_0)=0,~~~~~~~~~~~~~~~~~~~~~~~~~~~other~cases.~~~~~~~~~\\[2pt]
       \end{array}
    \right.
\end{align}
\end{lem}
\begin{lem}{\rm \cite{Wa3}}
\begin{align}
\Gamma_{st}^k(x_0)&=\left\{
       \begin{array}{c}
        \Gamma^n_{ii}(x_0)=\frac{1}{2}h'(0),~~~~~~~~~~ ~~~~~~~~~~~{\rm if }~s=t=i,k=n,i<n; \\[2pt]
        \Gamma^i_{ni}(x_0)=-\frac{1}{2}h'(0),~~~~~~~~~~~~~~~~~~~{\rm if }~s=n,t=i,k=i,i<n;\\[2pt]
        \Gamma^i_{in}(x_0)=-\frac{1}{2}h'(0),~~~~~~~~~~~~~~~~~~~{\rm if }~s=i,t=n,k=i,i<n,\\[2pt]
        \Gamma_{st}^i(x_0)=0,~~~~~~~~~~~~~~~~~~~~~~~~~~~other~cases.~~~~~~~~~
       \end{array}
    \right.
\end{align}
\end{lem}

When $n$ is even, then ${\rm tr}[{\rm \texttt{id}}]=2^{\frac{n}{2}},$ where tr as shorthand of ${\rm
trace}_{S(TM)}$, since the sum is taken over $r+l-k-j-|\alpha|-1=-n,~~r\leq -1,l\leq-n+3,$ then we have the following five cases:\\

\noindent {\bf case a)~I)}~$r=-1,~l=-n+3,~k=j=0,~|\alpha|=1$.\\

By applying the formula shown in (\ref{c2}), we can calculate
\begin{equation}
\Phi_1=-\int_{|\xi'|=1}\int^{+\infty}_{-\infty}\sum_{|\alpha|=1}{\rm tr}[\partial^\alpha_{\xi'}\pi^+_{\xi_n}\sigma_{-1}({D}_{J}^{-1})\times\partial^\alpha_{x'}\partial_{\xi_n}\sigma_{-n+3}({D}_{J}^{-n+3})](x_0)d\xi_n\sigma(\xi')dx'.
\end{equation}
For $i<n,$ (3.39) in \cite{LSW1} makes it obvious that
\begin{align}
\pi^+_{\xi_n}\partial_{\xi_i}\left(\frac{ic[J(\xi)]}{|\xi|^2}\right)(x_0)|_{|\xi'|=1}
&=\frac{1}{2(\xi_n-i)}c[J(dx_i)]+\frac{2i-\xi_n}{2(\xi_n-i)^2}\sum^{n-1}_{q=1}\xi_{i}\xi_{q}c[J(dx_q)]\\
&-\frac{1}{2(\xi_n-i)^2}\xi_{i}c[J(dx_n)]\nonumber\\
&=\frac{1}{2(\xi_n-i)}\sum^{n}_{\alpha=1}a^{i}_{\alpha}c(dx_{\alpha})+\frac{2i-\xi_n}{2(\xi_n-i)^2}\sum^{n}_{\beta=1}\sum^{n-1}_{q=1}\xi_{i}\xi_{q}a^{q}_{\beta}c(dx_{\beta})\nonumber\\
&-\frac{1}{2(\xi_n-i)^2}\sum^{n}_{\gamma=1}\xi_{i}a^{n}_{\gamma}c(dx_{\gamma}),\nonumber
\end{align}
where $c[J(dx_{p})]=\sum^{n}_{h=1}a^{p}_{h}c(dx_{h})$.

We get
\begin{align}
\partial_{x_i}\left(ic[J(\xi)]|\xi|^{-n+2}\right)(x_0)
=\frac{i\partial_{x_i}(c[J(\xi)])(x_0)}{|\xi|^{n-2}}-\frac{ic[J(\xi)]\partial_{x_i}(|\xi|^{n-2})(x_0)}{|\xi|^{2n-4}}
=\frac{i\partial_{x_i}(c[J(\xi)])(x_0)}{|\xi|^{n-2}},
\end{align}
of course,
\begin{align}
\partial_{x_i}\left(\frac{ic[J(\xi)]}{|\xi|^{n-2}}\right)(x_0)
&=\frac{i\partial_{x_i}(c[J(\sum^{n}_{p=1}\xi_{p}dx_{p})])(x_0)}{|\xi|^{n-2}}=\frac{i\sum^{n}_{p=1}\xi_{p}\partial_{x_i}(c(\sum^{n}_{h=1}a^{p}_{h}dx_{h}))(x_0)}{|\xi|^{n-2}}\\
&=\frac{i\sum^{n}_{p,h=1}\xi_{p}\partial_{x_i}(a^{p}_{h})c(dx_{h})(x_0)}{|\xi|^{n-2}}+\frac{i\sum^{n}_{p,h=1}\xi_{p}a^{p}_{h}\partial_{x_i}(c(dx_{h}))(x_0)}{|\xi|^{n-2}}\nonumber\\
&=\frac{i\sum^{n}_{p,h=1}\xi_{p}\partial_{x_i}(a^{p}_{h})c(dx_{h})(x_0)}{|\xi|^{n-2}}.\nonumber
\end{align}
When $|\xi'|=1$, we see that
\begin{align}
\partial_{x_i}\left(\frac{ic[J(\xi)]}{|\xi|^{n-2}}\right)(x_0)|_{|\xi'|=1}
&=\frac{i\sum^{n}_{h=1}\sum^{n-1}_{p=1}\xi_{p}\partial_{x_i}(a^{p}_{h})c(dx_{h})(x_0)}{(1+\xi_{n}^{2})^{\frac{n}{2}-1}}+\frac{i\sum^{n}_{h=1}\xi_{n}\partial_{x_i}(a^{n}_{h})c(dx_{h})(x_0)}{(1+\xi_{n}^{2})^{\frac{n}{2}-1}}.
\end{align}
An easy computation shows that
\begin{align}
\partial_{\xi_{n}}\partial_{x_i}\left(\frac{ic[J(\xi)]}{|\xi|^{n-2}}\right)(x_0)|_{|\xi'|=1}
&=i\sum^{n}_{h=1}\sum^{n-1}_{p=1}\xi_{p}\partial_{x_i}(a^{p}_{h})c(dx_{h})\partial_{\xi_{n}}\left(\frac{1}{(1+\xi_{n}^{2})^{\frac{n}{2}-1}}\right)\end{align}
\begin{align}
&+i\sum^{n}_{h=1}\partial_{x_i}(a^{n}_{h})c(dx_{h})\partial_{\xi_{n}}\left(\frac{\xi_{n}}{(1+\xi_{n}^{2})^{\frac{n}{2}-1}}\right)\nonumber\\
&=-\frac{i(n-2)\xi_{n}}{(1+\xi_{n}^{2})^{\frac{n}{2}}}\sum^{n}_{h=1}\sum^{n-1}_{p=1}\xi_{p}\partial_{x_i}(a^{p}_{h})c(dx_{h})\nonumber\\
&+\frac{i[1-(n-3)\xi_{n}^2]}{(1+\xi_{n}^{2})^{\frac{n}{2}}}\sum^{n}_{h=1}\partial_{x_i}(a^{n}_{h})c(dx_{h}).\nonumber
\end{align}

Hence, we have
\begin{align}
&\sum_{|\alpha|=1}{\rm tr}[\partial^\alpha_{\xi'}\pi^+_{\xi_n}\sigma_{-1}({D}_{J}^{-1})\times\partial^\alpha_{x'}\partial_{\xi_n}\sigma_{-n+3}({D}_{J}^{-n+3})](x_0)|_{|\xi'|=1}\\
=&-\frac{i(n-2)\xi _n}{2\left(\xi _n-i\right)^{\frac{n}{2}+1} \left(\xi _n+i\right)^{\frac{n}{2}}}\sum_{\alpha,h=1}^{n}\sum_{i,p=1}^{n-1}{\rm tr}[\xi_{p}a_{\alpha}^{i}\partial_{x_i}(a_{h}^{p})c(dx_{\alpha})c(dx_{h})]\nonumber\\
&+\frac{i[1-(n-3)\xi_n^2]}{2 \left(\xi _n-i\right)^{\frac{n}{2}+1} \left(\xi _n+i\right)^{\frac{n}{2}}}\sum_{\alpha,h=1}^{n}\sum_{i=1}^{n-1}{\rm tr}[a_{\alpha}^{i}\partial_{x_i}(a_{h}^{n})c(dx_{\alpha})c(dx_{h})]\nonumber\\
&-\frac{i(2i-\xi_{n})(n-2)\xi_n}{2\left(\xi_n-i\right)^{\frac{n}{2}+2} \left(\xi_n+i\right)^{\frac{n}{2}}}\sum_{\beta,h=1}^{n}\sum_{i,q,p=1}^{n-1}{\rm tr}[\xi_{i}\xi_{q}\xi_{p}a_{\beta}^{q}\partial_{x_i}(a_{h}^{p})c(dx_{\beta})c(dx_{h})]\nonumber\\
&+\frac{i(2i-\xi_{n})[1-(n-3)\xi_n^2]}{2 \left(\xi _n-i\right)^{\frac{n}{2}+2} \left(\xi _n+i\right)^{\frac{n}{2}}}\sum_{\beta,h=1}^{n}\sum_{i,q=1}^{n-1}{\rm tr}[\xi_{i}\xi_{q}a_{\beta}^{q}\partial_{x_i}(a_{h}^{n})c(dx_{\beta})c(dx_{h})]\nonumber\\
&+\frac{i(n-2)\xi_n}{2\left(\xi_n-i\right)^{\frac{n}{2}+2} \left(\xi _n+i\right)^{\frac{n}{2}}}\sum_{\gamma,h=1}^{n}\sum_{i,p=1}^{n-1}{\rm tr}[\xi_{i}\xi_{p}a_{\gamma}^{n}\partial_{x_i}(a_{h}^{p})c(dx_{\gamma})c(dx_{h})]\nonumber\\
&-\frac{i[1-(n-3)\xi_n^2]}{2 \left(\xi _n-i\right)^{\frac{n}{2}+2} \left(\xi _n+i\right)^{\frac{n}{2}}}\sum_{\gamma,h=1}^{n}\sum_{i=1}^{n-1}{\rm tr}[\xi_{i}a_{\gamma}^{n}\partial_{x_i}(a_{h}^{n})c(dx_{\gamma})c(dx_{h})].\nonumber
\end{align}

We note that $\int_{|\xi'|=1}{\{\xi_{i_1}\cdot\cdot\cdot\xi_{i_{2d+1}}}\}\sigma(\xi')=0,$ this gives
\begin{align}
\Phi_1&=-\int_{|\xi'|=1}\int^{+\infty}_{-\infty}\sum_{|\alpha|=1}{\rm tr}
[\partial^{\alpha}_{\xi'}\pi^{+}_{\xi_{n}}\sigma_{-1}({D}_{J}^{-1})\times\partial^{\alpha}_{x'}\partial_{\xi_{n}}\sigma_{-n+3}({D}_{J}^{-n+3})](x_0)d\xi_n\sigma(\xi')dx'\\
&=-\int_{|\xi'|=1}\int^{+\infty}_{-\infty}\frac{i[1-(n-3)\xi_n^2]}{2 \left(\xi _n-i\right)^{\frac{n}{2}+1} \left(\xi _n+i\right)^{\frac{n}{2}}}\sum_{\alpha,h=1}^{n}\sum_{i=1}^{n-1}{\rm tr}[a_{\alpha}^{i}\partial_{x_i}(a_{h}^{n})c(dx_{\alpha})c(dx_{h})]d\xi_n\sigma(\xi')dx'\nonumber\\
&-\int_{|\xi'|=1}\int^{+\infty}_{-\infty}\frac{i(2i-\xi_{n})[1-(n-3)\xi_n^2]}{2 \left(\xi _n-i\right)^{\frac{n}{2}+2} \left(\xi _n+i\right)^{\frac{n}{2}}}\sum_{\beta,h=1}^{n}\sum_{i,q=1}^{n-1}{\rm tr}[\xi_{i}\xi_{q}a_{\beta}^{q}\partial_{x_i}(a_{h}^{n})c(dx_{\beta})c(dx_{h})]d\xi_n\sigma(\xi')dx'\nonumber
\end{align}
\begin{align}
&-\int_{|\xi'|=1}\int^{+\infty}_{-\infty}\frac{i(n-2)\xi_n}{2\left(\xi_n-i\right)^{\frac{n}{2}+2} \left(\xi _n+i\right)^{\frac{n}{2}}}\sum_{\gamma,h=1}^{n}\sum_{i,p=1}^{n-1}{\rm tr}[\xi_{i}\xi_{p}a_{\gamma}^{n}\partial_{x_i}(a_{h}^{p})c(dx_{\gamma})c(dx_{h})]d\xi_n\sigma(\xi')dx'\nonumber.
\end{align}

Because $c(e_i)c(e_j)+c(e_j)c(e_i)=-2\delta_i^j$ then by the relation of the Clifford action and ${\rm tr}{AB}={\rm tr}{BA}$, we have the following equalities:
\begin{align}
&\sum_{\alpha,h=1}^{n}\sum_{i=1}^{n-1}{\rm tr}[a_{\alpha}^{i}\partial_{x_i}(a_{h}^{n})c(dx_{\alpha})c(dx_{h})]=-\sum_{h=1}^{n}\sum_{i=1}^{n-1}a_{h}^{i}\partial_{x_i}(a_{h}^{n}){\rm tr}[\texttt{id}];\\
&\sum_{\beta,h=1}^{n}\sum_{i,q=1}^{n-1}{\rm tr}[\xi_{i}\xi_{q}a_{\beta}^{q}\partial_{x_i}(a_{h}^{n})c(dx_{\beta})c(dx_{h})]=-\sum_{h=1}^{n}\sum_{i,q=1}^{n-1}\xi_{i}\xi_{q}a_{h}^{q}\partial_{x_i}(a_{h}^{n}){\rm tr}[\texttt{id}];\\
&\sum_{\gamma,h=1}^{n}\sum_{i,p=1}^{n-1}{\rm tr}[\xi_{i}\xi_{p}a_{\gamma}^{n}\partial_{x_i}(a_{h}^{p})c(dx_{\gamma})c(dx_{h})]=-\sum_{h=1}^{n}\sum_{i,p=1}^{n-1}\xi_{i}\xi_{p}a_{h}^{n}\partial_{x_i}(a_{h}^{p}){\rm tr}[\texttt{id}].
\end{align}
And by $\int_{|\xi'|=1}\sum_{i,j=1}^{n-1}\xi_i\xi_j\sigma(\xi')=\frac{1}{n-1}\sum_{i,j=1}^{n-1}\delta_{j}^{i}Vol(S^{n-2})$ in \cite{LW}, we have
\begin{align}
\Phi_1
&=\int^{+\infty}_{-\infty}\frac{i[1-(n-3)\xi_n^2]}{2 \left(\xi _n-i\right)^{\frac{n}{2}+1} \left(\xi _n+i\right)^{\frac{n}{2}}}\sum_{h=1}^{n}\sum_{i=1}^{n-1}a_{h}^{i}\partial_{x_i}(a_{h}^{n})\cdot Vol(S^{n-2}){\rm tr}[\texttt{id}]d\xi_ndx'\\
&+\int^{+\infty}_{-\infty}\frac{i(2i-\xi_{n})[1-(n-3)\xi_n^2]}{2 \left(\xi _n-i\right)^{\frac{n}{2}+2} \left(\xi _n+i\right)^{\frac{n}{2}}}\sum_{h=1}^{n}\sum_{i=1}^{n-1}a_{h}^{i}\partial_{x_i}(a_{h}^{n})\cdot \frac{1}{n-1}Vol(S^{n-2}){\rm tr}[\texttt{id}]d\xi_ndx'\nonumber\\
&+\int^{+\infty}_{-\infty}\frac{i(n-2)\xi_n}{2\left(\xi_n-i\right)^{\frac{n}{2}+2} \left(\xi _n+i\right)^{\frac{n}{2}}}\sum_{h=1}^{n}\sum_{i=1}^{n-1}a_{h}^{n}\partial_{x_i}(a_{h}^{i})\cdot \frac{1}{n-1}Vol(S^{n-2}){\rm tr}[\texttt{id}]d\xi_ndx'.\nonumber
\end{align}

It is easily seen that
\begin{align}
\Phi_1
&=\sum_{h=1}^{n}\sum_{i=1}^{n-1}a_{h}^{i}\partial_{x_i}(a_{h}^{n})\cdot Vol(S^{n-2}){\rm tr}[\texttt{id}]\cdot \frac{2\pi i}{(\frac{n}{2})!}\mathcal{A}_{0}dx'\\
&+\sum_{h=1}^{n}\sum_{i=1}^{n-1}a_{h}^{i}\partial_{x_i}(a_{h}^{n})\cdot \frac{1}{n-1}Vol(S^{n-2}){\rm tr}[\texttt{id}]\cdot \frac{2\pi i}{(\frac{n}{2}+1)!}\mathcal{A}_{1}dx'\nonumber\\
&+\sum_{h=1}^{n}\sum_{i=1}^{n-1}a_{h}^{n}\partial_{x_i}(a_{h}^{i})\cdot \frac{1}{n-1}Vol(S^{n-2}){\rm tr}[\texttt{id}]\cdot \frac{2\pi i}{(\frac{n}{2}+1)!}\mathcal{A}_{2}dx',\nonumber
\end{align}
where
\begin{align}
\mathcal{A}_{0}
&=\Big[\frac{-i(n-3)\xi_n^2+i}{2\left(\xi _n+i\right)^{\frac{n}{2}}}\Big]^{(\frac{n}{2})}\bigg|_{\xi_n=i};\\
\mathcal{A}_{1}
&=\Big[\frac{i(n-3)\xi_n^3+2(n-3)\xi_n^2-i\xi_n-2}{2\left(\xi _n+i\right)^{\frac{n}{2}}}\Big]^{(\frac{n}{2}+1)}\bigg|_{\xi_n=i};
\end{align}
\begin{align}
\mathcal{A}_{2}
&=\Big[\frac{i(n-2)\xi_n}{2\left(\xi _n+i\right)^{\frac{n}{2}}}\Big]^{(\frac{n}{2}+1)}\bigg|_{\xi_n=i}.
\end{align}

\noindent {\bf case (a)~(II)}~$r=-1, l=-n+3, |\alpha|=k=0, j=1$.\\

It is easy to check that
\begin{equation}
\Phi_2=-\frac{1}{2}\int_{|\xi'|=1}\int^{+\infty}_{-\infty}{\rm tr}[\partial_{x_{n}}\pi^{+}_{\xi_{n}}\sigma_{-1}({D}_{J}^{-1})\times\partial^{2}_{\xi_{n}}\sigma_{-n+3}({D}_{J}^{-n+3})](x_0)d\xi_n\sigma(\xi')dx'.
\end{equation}
We can assert that
\begin{align}
&\pi^+_{\xi_n}\partial_{x_n}\left(\frac{ic[J(\xi)]}{|\xi|^2}\right)(x_0)|_{|\xi'|=1}\\
&=\frac{1}{2(\xi_n-i)}\sum^{n}_{h=1}\sum^{n-1}_{p=1}\xi_{p}\partial_{x_n}(a^{p}_{h})c(dx_{h})+\frac{i}{2(\xi_n-i)}\sum^{n}_{h=1}\partial_{x_n}(a^{n}_{h})c(dx_{h})\nonumber\\
&+\frac{1}{2(\xi_n-i)}\sum^{n-1}_{h,p=1}\xi_{p}a^{p}_{h}\partial_{x_n}(c(dx_{h}))+\frac{i}{2(\xi_n-i)}\sum^{n-1}_{h=1}a^{n}_{h}\partial_{x_n}(c(dx_{h}))\nonumber\\
&+\frac{2i-\xi_n}{4(\xi_n-i)^2}h'(0)\sum^{n}_{h=1}\sum^{n-1}_{p=1}\xi_{p}a^{p}_{h}c(dx_{h})-\frac{1}{4(\xi_n-i)^2}h'(0)\sum^{n}_{h=1}a^{n}_{h}c(dx_{h}),\nonumber
\end{align}
where $\sum_{h=1}^{n-1}\partial_{x_n}(c(dx_h))=\sum_{h=1}^{n-1}\frac{1}{2}h'(0)c(dx_h).$

By calculation, we have
\begin{align}
\partial_{\xi_n}^2\left(\frac{ic[J(\xi)]}{|\xi|^{n-2}}\right)(x_0)
&=\partial_{\xi_n}^2\left(\frac{i\sum^{n}_{i=1}\xi_{i}c[J(dx_i)]}{|\xi|^{n-2}}\right)(x_0)\\
&=\partial_{\xi_n}^2\left(\frac{i\sum^{n}_{i,\beta=1}\xi_{i}a_{\beta}^ic(dx_{\beta})}{|\xi|^{n-2}}\right)(x_0)\nonumber
\end{align}
and
\begin{align}
\partial_{\xi_n}^2\left(\frac{ic[J(\xi)]}{|\xi|^{n-2}}\right)(x_0)|_{|\xi'|=1}
&=i\sum^{n}_{\beta=1}\sum^{n-1}_{i=1}\xi_{i}a_{\beta}^ic(dx_{\beta})\partial_{\xi_n}^2\left(\frac{1}{(1+\xi_{n}^{2})^{\frac{n}{2}-1}}\right)\\
&+i\sum^{n}_{\beta=1}a_{\beta}^nc(dx_{\beta})\partial_{\xi_n}^2\left(\frac{\xi_{n}}{(1+\xi_{n}^{2})^{\frac{n}{2}-1}}\right)\nonumber\\
&=\frac{i(n-2)[-1+(n-1)\xi_n^2]}{\left(1+\xi _n^2\right)^{\frac{n}{2}+1}}\sum^{n}_{\beta=1}\sum^{n-1}_{i=1}\xi_{i}a_{\beta}^ic(dx_{\beta})\nonumber\\
&+\frac{i(n-2)\xi_n[-3+(n-3)\xi_n^2]}{\left(1+\xi _n^2\right)^{\frac{n}{2}+1}}\sum^{n}_{\beta=1}a_{\beta}^nc(dx_{\beta}).\nonumber
\end{align}

A trivial verification shows that
\begin{align}\label{c3}
&\Phi_2
=-\frac{1}{2}\int_{|\xi'|=1}\int^{+\infty}_{-\infty}{\rm tr}[\partial_{x_{n}}\pi^{+}_{\xi_{n}}\sigma_{-1}({D}_{J}^{-1})\times\partial^{2}_{\xi_{n}}\sigma_{-n+3}({D}_{J}^{-n+3})](x_0)d\xi_n\sigma(\xi')dx'\\
&=-\frac{1}{2}\int_{|\xi'|=1}\int^{+\infty}_{-\infty}\frac{i(n-2)[-1+(n-1)\xi_n^{2}]}{2 \left(\xi _n-i\right)^{\frac{n}{2}+2} \left(\xi _n+i\right)^{\frac{n}{2}+1}}\sum_{h,\beta=1}^{n}\sum_{p,i=1}^{n-1}{\rm tr}[\xi_{p}\xi_{i}a_{\beta}^{i}\partial_{x_n}(a_{h}^{p})c(dx_{h})c(dx_{\beta})]d\xi_n\sigma(\xi')dx'\nonumber\\
&-\frac{1}{2}\int_{|\xi'|=1}\int^{+\infty}_{-\infty}-\frac{(n-2)\xi_n[-3+(n-3)\xi_n^{2}]}{2 \left(\xi _n-i\right)^{\frac{n}{2}+2} \left(\xi _n+i\right)^{\frac{n}{2}+1}}\sum_{h,\beta=1}^{n}{\rm tr}[a_{\beta}^{n}\partial_{x_n}(a_{h}^{n})c(dx_{h})c(dx_{\beta})]d\xi_n\sigma(\xi')dx'\nonumber\\
&-\frac{1}{2}\int_{|\xi'|=1}\int^{+\infty}_{-\infty}\frac{i(n-2)[-1+(n-1)\xi_n^{2}]}{2 \left(\xi _n-i\right)^{\frac{n}{2}+2} \left(\xi _n+i\right)^{\frac{n}{2}+1}}\sum_{\beta=1}^{n}\sum_{h,p,i=1}^{n-1}{\rm tr}[\xi_{p}\xi_{i}a_{h}^{p}a_{\beta}^{i}\partial_{x_n}(c(dx_{h}))c(dx_{\beta})]d\xi_n\sigma(\xi')dx'\nonumber\\
&-\frac{1}{2}\int_{|\xi'|=1}\int^{+\infty}_{-\infty}-\frac{(n-2)\xi_n[-3+(n-3)\xi_n^{2}]}{2 \left(\xi _n-i\right)^{\frac{n}{2}+2} \left(\xi _n+i\right)^{\frac{n}{2}+1}}\sum_{\beta=1}^{n}\sum_{h=1}^{n-1}{\rm tr}[a_{h}^{n}a_{\beta}^{n}\partial_{x_n}(c(dx_{h}))c(dx_{\beta})]d\xi_n\sigma(\xi')dx'\nonumber\\
&-\frac{1}{2}\int_{|\xi'|=1}\int^{+\infty}_{-\infty}\frac{i(n-2)(2i-\xi_n)[-1+(n-1)\xi_n^{2}]}{4 \left(\xi _n-i\right)^{\frac{n}{2}+3} \left(\xi _n+i\right)^{\frac{n}{2}+1}}h'(0)\sum_{h,\beta=1}^{n}\sum_{p,i=1}^{n-1}{\rm tr}[\xi_{p}\xi_{i}a_{h}^{p}a_{\beta}^{i}c(dx_{h})c(dx_{\beta})]d\xi_n\sigma(\xi')dx'\nonumber\\
&-\frac{1}{2}\int_{|\xi'|=1}\int^{+\infty}_{-\infty}-\frac{i(n-2)\xi_n[-3+(n-3)\xi_n^{2}]}{4 \left(\xi _n-i\right)^{\frac{n}{2}+3} \left(\xi _n+i\right)^{\frac{n}{2}+1}}h'(0)\sum_{h,\beta=1}^{n}{\rm tr}[a_{h}^{n}a_{\beta}^{n}c(dx_{h})c(dx_{\beta})]d\xi_n\sigma(\xi')dx'.\nonumber
\end{align}
As in the proof of  case (a)~(I), equation (\ref{c3}) gives
\begin{align}
&\Phi_2
=-\frac{1}{2}\int^{+\infty}_{-\infty}-\frac{i(n-2)[-1+(n-1)\xi_n^{2}]}{2 \left(\xi _n-i\right)^{\frac{n}{2}+2} \left(\xi _n+i\right)^{\frac{n}{2}+1}}\sum_{h=1}^{n}\sum_{i=1}^{n-1}a_{h}^{i}\partial_{x_n}(a_{h}^{i})\cdot \frac{1}{n-1}Vol(S^{n-2}){\rm tr}[\texttt{id}]d\xi_ndx'\\
&-\frac{1}{2}\int^{+\infty}_{-\infty}\frac{(n-2)\xi_n[-3+(n-3)\xi_n^{2}]}{2 \left(\xi _n-i\right)^{\frac{n}{2}+2} \left(\xi _n+i\right)^{\frac{n}{2}+1}}\sum_{h=1}^{n}a_{h}^{n}\partial_{x_n}(a_{h}^{n})\cdot Vol(S^{n-2}){\rm tr}[\texttt{id}]d\xi_ndx'\nonumber\\
&-\frac{1}{2}\int^{+\infty}_{-\infty}-\frac{i(n-2)[-1+(n-1)\xi_n^{2}]}{4 \left(\xi _n-i\right)^{\frac{n}{2}+2} \left(\xi _n+i\right)^{\frac{n}{2}+1}}h'(0)\sum_{h,i=1}^{n-1}(a_{h}^{i})^{2}\cdot \frac{1}{n-1}Vol(S^{n-2}){\rm tr}[\texttt{id}]d\xi_ndx'\nonumber\\
&-\frac{1}{2}\int^{+\infty}_{-\infty}\frac{(n-2)\xi_n[-3+(n-3)\xi_n^{2}]}{4 \left(\xi _n-i\right)^{\frac{n}{2}+2} \left(\xi _n+i\right)^{\frac{n}{2}+1}}h'(0)\sum_{h=1}^{n-1}(a_{h}^{n})^{2}\cdot Vol(S^{n-2}){\rm tr}[\texttt{id}]d\xi_ndx'\nonumber\\
&-\frac{1}{2}\int^{+\infty}_{-\infty}-\frac{i(n-2)(2i-\xi_n)[-1+(n-1)\xi_n^{2}]}{4 \left(\xi _n-i\right)^{\frac{n}{2}+3} \left(\xi _n+i\right)^{\frac{n}{2}+1}}h'(0)\sum_{h=1}^{n}\sum_{i=1}^{n-1}(a_{h}^{i})^{2}\cdot \frac{1}{n-1}Vol(S^{n-2}){\rm tr}[\texttt{id}]d\xi_ndx'\nonumber
\end{align}
\begin{align}
&-\frac{1}{2}\int^{+\infty}_{-\infty}\frac{i(n-2)\xi_n[-3+(n-3)\xi_n^{2}]}{4 \left(\xi _n-i\right)^{\frac{n}{2}+3} \left(\xi _n+i\right)^{\frac{n}{2}+1}}h'(0)\sum_{h=1}^{n}(a_{h}^{n})^{2}\cdot Vol(S^{n-2}){\rm tr}[\texttt{id}]d\xi_ndx'.\nonumber
\end{align}

We can certainly assume that $J(\frac{\partial}{\partial{x_{l}}})=\sum^{n}_{q=1}\widetilde{a}^{q}_{l}\frac{\partial}{\partial{x_{q}}},$ since
\begin{align}
\langle \frac{\partial}{\partial{x_{l}}}, J(dx_{p})\rangle&=\langle J(\frac{\partial}{\partial{x_{l}}}), dx_{p}\rangle,
\end{align}
then we have $a^{q}_{h}=\widetilde{a}^{q}_{h}.$
As
\begin{align}
J^2(x_0)=\texttt{id}, J^{T}J(x_0)=\texttt{id},
\end{align}
we have $a^{p}_{l}a^{j}_{p}(x_0)=\delta_l^j, a^{p}_{l}(x_0)=a^{l}_{p}(x_0)$.
Let $J(e_l)=\sum^{n}_{q=1}{b}^{q}_{l}e_q$, then $\partial{x_{j}}({b}^{q}_{l})(x_0)=\partial{x_{j}}({a}^{q}_{l})(x_0)$ and ${b}^{q}_{l}(x_0)={a}^{q}_{l}(x_0)$.
We next show that
\begin{align}
0&=\sum_{i=1}^{n-1}\partial_{x_n}(\delta_{i}^{i})=\sum_{h=1}^{n}\sum_{i=1}^{n-1}\partial_{x_n}(a_{h}^{i}a_{h}^{i})=2\sum_{h=1}^{n}\sum_{i=1}^{n-1}a_{h}^{i}\partial_{x_n}(a_{h}^{i});\\
0&=\partial_{x_n}(\delta_{n}^{n})=\sum_{h=1}^{n}\partial_{x_n}(a_{h}^{n}a_{h}^{n})=2\sum_{h=1}^{n}a_{h}^{n}\partial_{x_n}(a_{h}^{n}).
\end{align}

Therefore
\begin{align}
&\Phi_2
=h'(0)\sum_{h,i=1}^{n-1}(a_{h}^{i})^{2}\cdot \frac{1}{n-1}Vol(S^{n-2}){\rm tr}[\texttt{id}]\cdot \frac{2\pi i}{(\frac{n}{2}+1)!}\mathcal{A}_{3}dx'\\
&+h'(0)\sum_{i=1}^{n-1}(a_{n}^{i})^{2}\cdot Vol(S^{n-2}){\rm tr}[\texttt{id}]\cdot \frac{2\pi i}{(\frac{n}{2}+1)!}\mathcal{A}_{4}dx'\nonumber\\
&+h'(0)\sum_{h=1}^{n}\sum_{i=1}^{n-1}(a_{h}^{i})^{2}\cdot \frac{1}{n-1}Vol(S^{n-2}){\rm tr}[\texttt{id}]\cdot \frac{2\pi i}{(\frac{n}{2}+2)!}\mathcal{A}_{5}dx'\nonumber\\
&+h'(0)\sum_{h=1}^{n}(a_{h}^{n})^{2}\cdot Vol(S^{n-2}){\rm tr}[\texttt{id}]\cdot \frac{2\pi i}{(\frac{n}{2}+2)!}\mathcal{A}_{6}dx',\nonumber
\end{align}
where
\begin{align}
\mathcal{A}_{3}
&=\Big[\frac{i(n^{2}-3n+2)\xi_n^{2}-i(n-2)}{8\left(\xi _n+i\right)^{\frac{n}{2}+1}}\Big]^{(\frac{n}{2}+1)}\bigg|_{\xi_n=i};\\
\mathcal{A}_{4}
&=\Big[\frac{-(n^{2}-5n+6)\xi_n^{3}+3(n-2)\xi_n}{8\left(\xi _n+i\right)^{\frac{n}{2}+1}}\Big]^{(\frac{n}{2}+1)}\bigg|_{\xi_n=i};\\
\mathcal{A}_{5}
&=\Big[\frac{-i(n^{2}-3n+2)\xi_n^{3}-2(n^{2}-3n+2)\xi_n^{2}+i(n-2)\xi_n+2(n-2)}{8\left(\xi _n+i\right)^{\frac{n}{2}+1}}\Big]^{(\frac{n}{2}+2)}\bigg|_{\xi_n=i};\\
\mathcal{A}_{6}
&=\Big[\frac{-i(n^{2}-5n+6)\xi_n^{3}+3i(n-2)\xi_n}{8\left(\xi _n+i\right)^{\frac{n}{2}+1}}\Big]^{(\frac{n}{2}+2)}\bigg|_{\xi_n=i}.
\end{align}

\noindent  {\bf case (a)~(III)}~$r=-1,l=-n+3,|\alpha|=j=0,k=1$.\\

Using (\ref{c2}), we get
\begin{equation}
\Phi_3=-\frac{1}{2}\int_{|\xi'|=1}\int^{+\infty}_{-\infty}{\rm tr}[\partial_{\xi_{n}}\pi^{+}_{\xi_{n}}\sigma_{-1}({D}_{J}^{-1})
      \times\partial_{\xi_{n}}\partial_{x_{n}}\sigma_{-n+3}({D}_{J}^{-n+3})](x_0)d\xi_n\sigma(\xi')dx'.
\end{equation}\\

We check at once that
\begin{align}
\pi^+_{\xi_n}\partial_{\xi_n}\left(\frac{ic[J(\xi)]}{|\xi|^2}\right)(x_0)|_{|\xi'|=1}
&=-\frac{1}{2(\xi_{n}-i)^2}\sum^{n}_{\beta=1}\sum^{n-1}_{i=1}\xi_{i}a_{\beta}^ic(dx_{\beta})-\frac{i}{2(\xi_{n}-i)^2}\sum^{n}_{\beta=1}a_{\beta}^nc(dx_{\beta}).
\end{align}

Likewise,
\begin{align}
\partial_{x_n}\left(\frac{ic[J(\xi)]}{|\xi|^{n-2}}\right)(x_0)
&=\frac{i\partial_{x_n}(c[J(\xi)])(x_0)}{|\xi|^{n-2}}-\frac{ic[J(\xi)]\partial_{x_n}(|\xi|^{n-2})(x_0)}{|\xi|^{2n-4}}\\
&=\frac{i\sum^{n}_{p,h=1}\xi_{p}\partial_{x_n}(a^{p}_{h})c(dx_{h})(x_0)}{|\xi|^{n-2}}+\frac{i\sum^{n}_{p,h=1}\xi_{p}a^{p}_{h}\partial_{x_n}(c(dx_{h}))(x_0)}{|\xi|^{n-2}}\nonumber\\
&-\frac{i(\frac{n}{2}-1)h'(0)|\xi^{'}|^{2}\sum^{n}_{p,h=1}\xi_{p}a^{p}_{h}c(dx_{h})(x_0)}{|\xi|^{n}}.\nonumber
\end{align}
Then we have
\begin{align}
\partial_{x_n}\left(\frac{ic[J(\xi)]}{|\xi|^{n-2}}\right)(x_0)|_{|\xi'|=1}
&=\frac{i\sum^{n}_{p,h=1}\xi_{p}\partial_{x_n}(a^{p}_{h})c(dx_{h})(x_0)}{(1+\xi_{n}^{2})^{\frac{n}{2}-1}}+\frac{i\sum^{n}_{p,h=1}\xi_{p}a^{p}_{h}\partial_{x_n}(c(dx_{h}))(x_0)}{(1+\xi_{n}^{2})^{\frac{n}{2}-1}}\\
&-\frac{i(\frac{n}{2}-1)h'(0)\sum^{n}_{p,h=1}\xi_{p}a^{p}_{h}c(dx_{h})(x_0)}{(1+\xi_{n}^{2})^{\frac{n}{2}}}.\nonumber
\end{align}
Thus
\begin{align}
&\partial_{\xi_n}\partial_{x_n}\left(\frac{ic[J(\xi)]}{|\xi|^{n-2}}\right)(x_0)|_{|\xi'|=1}=-\frac{i(n-2)\xi _n}{\left(1+\xi _n^2\right)^{\frac{n}{2}}}\sum^{n}_{h=1}\sum^{n-1}_{p=1}\xi_{p}\partial_{x_n}(a^{p}_{h})c(dx_{h})\\
&+\frac{i[1-(n-3)\xi _n^2]}{\left(1+\xi _n^2\right)^{\frac{n}{2}}}\sum^{n}_{h=1}\partial_{x_n}(a^{n}_{h})c(dx_{h})-\frac{i(n-2)\xi _n}{\left(1+\xi _n^2\right)^{\frac{n}{2}}}\sum^{n-1}_{h,p=1}\xi_{p}a^{p}_{h}\partial_{x_n}(c(dx_{h}))\nonumber\\
&+\frac{i[1-(n-3)\xi _n^2]}{\left(1+\xi _n^2\right)^{\frac{n}{2}}}\sum^{n-1}_{h=1}a^{n}_{h}\partial_{x_n}(c(dx_{h}))+\frac{i(\frac{n}{2}-1)n\xi _n}{\left(1+\xi _n^2\right)^{\frac{n}{2}+1}}h'(0)\sum^{n}_{h=1}\sum^{n-1}_{p=1}\xi_{p}a^{p}_{h}c(dx_{h})\nonumber\\
&-\frac{i(\frac{n}{2}-1)[1-(n-1)\xi _n^2]}{\left(1+\xi _n^2\right)^{\frac{n}{2}+1}}h'(0)\sum^{n}_{h=1}a^{n}_{h}c(dx_{h}).\nonumber
\end{align}

We see at once that
\begin{align}
&\Phi_3=-\frac{1}{2}\int_{|\xi'|=1}\int^{+\infty}_{-\infty}{\rm tr}[\partial_{\xi_{n}}\pi^{+}_{\xi_{n}}\sigma_{-1}({D}_{J}^{-1})
      \times\partial_{\xi_{n}}\partial_{x_{n}}\sigma_{-n+3}({D}_{J}^{-n+3})](x_0)d\xi_n\sigma(\xi')dx'\\
&=-\frac{1}{2}\int_{|\xi'|=1}\int^{+\infty}_{-\infty}\frac{i(n-2)\xi_n}{2 \left(\xi _n-i\right)^{\frac{n}{2}+2} \left(\xi _n+i\right)^{\frac{n}{2}}}\sum_{h,\beta=1}^{n}\sum_{p,i=1}^{n-1}{\rm tr}[\xi_{p}\xi_{i}a_{\beta}^{i}\partial_{x_n}(a_{h}^{p})c(dx_{h})c(dx_{\beta})]d\xi_n\sigma(\xi')dx'\nonumber\\
&-\frac{1}{2}\int_{|\xi'|=1}\int^{+\infty}_{-\infty}\frac{i(n-2)\xi_n}{2 \left(\xi _n-i\right)^{\frac{n}{2}+2} \left(\xi _n+i\right)^{\frac{n}{2}}}\sum_{\beta=1}^{n}\sum_{h,p,i=1}^{n-1}{\rm tr}[\xi_{p}\xi_{i}a_{h}^{p}a_{\beta}^{i}\partial_{x_n}(c(dx_{h}))c(dx_{\beta})]d\xi_n\sigma(\xi')dx'\nonumber\\
&-\frac{1}{2}\int_{|\xi'|=1}\int^{+\infty}_{-\infty}-\frac{i(\frac{n}{2}-1)n\xi _n}{2 \left(\xi _n-i\right)^{\frac{n}{2}+3} \left(\xi _n+i\right)^{\frac{n}{2}+1}}h'(0)\sum_{h,\beta=1}^{n}\sum_{p,i=1}^{n-1}{\rm tr}[\xi_{p}\xi_{i}a_{h}^{p}a_{\beta}^{i}c(dx_{h})c(dx_{\beta})]d\xi_n\sigma(\xi')dx'\nonumber\\
&-\frac{1}{2}\int_{|\xi'|=1}\int^{+\infty}_{-\infty}\frac{1-(n-3)\xi_n^{2}}{2 \left(\xi _n-i\right)^{\frac{n}{2}+2} \left(\xi _n+i\right)^{\frac{n}{2}}}\sum_{h,\beta=1}^{n}{\rm tr}[a_{\beta}^{n}\partial_{x_n}(a_{h}^{n})c(dx_{h})c(dx_{\beta})]d\xi_n\sigma(\xi')dx'\nonumber\\
&-\frac{1}{2}\int_{|\xi'|=1}\int^{+\infty}_{-\infty}\frac{1-(n-3)\xi_n^{2}}{2 \left(\xi _n-i\right)^{\frac{n}{2}+2} \left(\xi _n+i\right)^{\frac{n}{2}}}\sum_{\beta=1}^{n}\sum_{h=1}^{n-1}{\rm tr}[a_{h}^{n}a_{\beta}^{n}\partial_{x_n}(c(dx_{h}))c(dx_{\beta})]d\xi_n\sigma(\xi')dx'\nonumber\\
&-\frac{1}{2}\int_{|\xi'|=1}\int^{+\infty}_{-\infty}-\frac{(\frac{n}{2}-1)[1-(n-1)\xi_n^{2}]}{2 \left(\xi _n-i\right)^{\frac{n}{2}+3} \left(\xi _n+i\right)^{\frac{n}{2}+1}}h'(0)\sum_{h,\beta=1}^{n}{\rm tr}[a_{h}^{n}a_{\beta}^{n}c(dx_{h})c(dx_{\beta})]d\xi_n\sigma(\xi')dx'.\nonumber
\end{align}
From this, we have
\begin{align}
&\Phi_3
=-\frac{1}{2}\int^{+\infty}_{-\infty}-\frac{i(n-2)\xi_n}{2 \left(\xi _n-i\right)^{\frac{n}{2}+2} \left(\xi _n+i\right)^{\frac{n}{2}}}\sum_{h=1}^{n}\sum_{i=1}^{n-1}a_{h}^{i}\partial_{x_n}(a_{h}^{i})\cdot \frac{1}{n-1}Vol(S^{n-2}){\rm tr}[\texttt{id}]d\xi_ndx'\\
&-\frac{1}{2}\int^{+\infty}_{-\infty}-\frac{i(n-2)\xi_n}{4 \left(\xi _n-i\right)^{\frac{n}{2}+2} \left(\xi _n+i\right)^{\frac{n}{2}}}h'(0)\sum_{h,i=1}^{n-1}(a_{h}^{i})^{2}\cdot \frac{1}{n-1}Vol(S^{n-2}){\rm tr}[\texttt{id}]d\xi_ndx'\nonumber\\
&-\frac{1}{2}\int^{+\infty}_{-\infty}\frac{i(\frac{n}{2}-1)n\xi _n}{2 \left(\xi _n-i\right)^{\frac{n}{2}+3} \left(\xi _n+i\right)^{\frac{n}{2}+1}}h'(0)\sum_{h=1}^{n}\sum_{i=1}^{n-1}(a_{h}^{i})^{2}\cdot \frac{1}{n-1}Vol(S^{n-2}){\rm tr}[\texttt{id}]d\xi_ndx'\nonumber\\
&-\frac{1}{2}\int^{+\infty}_{-\infty}-\frac{1-(n-3)\xi_n^{2}}{2 \left(\xi _n-i\right)^{\frac{n}{2}+2} \left(\xi _n+i\right)^{\frac{n}{2}}}\sum_{h=1}^{n}a_{h}^{n}\partial_{x_n}(a_{h}^{n})\cdot Vol(S^{n-2}){\rm tr}[\texttt{id}]d\xi_ndx'\nonumber\\
&-\frac{1}{2}\int^{+\infty}_{-\infty}-\frac{1-(n-3)\xi_n^{2}}{4 \left(\xi _n-i\right)^{\frac{n}{2}+2} \left(\xi _n+i\right)^{\frac{n}{2}}}h'(0)\sum_{h=1}^{n-1}(a_{h}^{n})^{2}\cdot Vol(S^{n-2}){\rm tr}[\texttt{id}]d\xi_ndx'\nonumber\\
&-\frac{1}{2}\int^{+\infty}_{-\infty}\frac{(\frac{n}{2}-1)[1-(n-1)\xi_n^{2}]}{2 \left(\xi _n-i\right)^{\frac{n}{2}+3} \left(\xi _n+i\right)^{\frac{n}{2}+1}}h'(0)\sum_{h=1}^{n}(a_{h}^{n})^{2}\cdot Vol(S^{n-2}){\rm tr}[\texttt{id}]d\xi_ndx'.\nonumber
\end{align}
It is immediate that
\begin{align}
&\Phi_3
=h'(0)\sum_{h,i=1}^{n-1}(a_{h}^{i})^{2}\cdot \frac{1}{n-1}Vol(S^{n-2}){\rm tr}[\texttt{id}]\cdot \frac{2\pi i}{(\frac{n}{2}+1)!}\mathcal{A}_{7}dx'\\
&+h'(0)\sum_{h=1}^{n}\sum_{i=1}^{n-1}(a_{h}^{i})^{2}\cdot \frac{1}{n-1}Vol(S^{n-2}){\rm tr}[\texttt{id}]\cdot \frac{2\pi i}{(\frac{n}{2}+2)!}\mathcal{A}_{8}dx'\nonumber\\
&+h'(0)\sum_{i=1}^{n-1}(a_{n}^{i})^{2}\cdot Vol(S^{n-2}){\rm tr}[\texttt{id}]\cdot \frac{2\pi i}{(\frac{n}{2}+1)!}\mathcal{A}_{9}dx'\nonumber\\
&+h'(0)\sum_{h=1}^{n}(a_{h}^{n})^{2}\cdot Vol(S^{n-2}){\rm tr}[\texttt{id}]\cdot \frac{2\pi i}{(\frac{n}{2}+2)!}\mathcal{A}_{10}dx',\nonumber
\end{align}
where
\begin{align}
\mathcal{A}_{7}
&=\Big[\frac{i(n-2)\xi_n}{8\left(\xi _n+i\right)^{\frac{n}{2}}}\Big]^{(\frac{n}{2}+1)}\bigg|_{\xi_n=i};\\
\mathcal{A}_{8}
&=\Big[\frac{-in(n-2)\xi_n}{8\left(\xi _n+i\right)^{\frac{n}{2}+1}}\Big]^{(\frac{n}{2}+2)}\bigg|_{\xi_n=i};\\
\mathcal{A}_{9}
&=\Big[\frac{-(n-3)\xi_n^{2}+1}{8\left(\xi _n+i\right)^{\frac{n}{2}}}\Big]^{(\frac{n}{2}+1)}\bigg|_{\xi_n=i};\\
\mathcal{A}_{10}
&=\Big[\frac{(n^{2}-3n+2)\xi_n^{2}-(n-2)}{8\left(\xi _n+i\right)^{\frac{n}{2}+1}}\Big]^{(\frac{n}{2}+2)}\bigg|_{\xi_n=i}.
\end{align}

In combination with the calculation,
\begin{align}
&\Phi_1+\Phi_2+\Phi_3\\
&=\sum_{h=1}^{n}\sum_{i=1}^{n-1}a_{h}^{i}\partial_{x_i}(a_{h}^{n})\cdot Vol(S^{n-2}){\rm tr}[\texttt{id}]\cdot \frac{2\pi i}{(\frac{n}{2})!}\mathcal{A}_{0}dx'\nonumber\\
&+\sum_{h=1}^{n}\sum_{i=1}^{n-1}a_{h}^{i}\partial_{x_i}(a_{h}^{n})\cdot \frac{1}{n-1}Vol(S^{n-2}){\rm tr}[\texttt{id}]\cdot \frac{2\pi i}{(\frac{n}{2}+1)!}\mathcal{A}_{1}dx'\nonumber\\
&+\sum_{h=1}^{n}\sum_{i=1}^{n-1}a_{h}^{n}\partial_{x_i}(a_{h}^{i})\cdot \frac{1}{n-1}Vol(S^{n-2}){\rm tr}[\texttt{id}]\cdot \frac{2\pi i}{(\frac{n}{2}+1)!}\mathcal{A}_{2}dx'\nonumber\\
&+h'(0)\sum_{h,i=1}^{n-1}(a_{h}^{i})^{2}\cdot \frac{1}{n-1}Vol(S^{n-2}){\rm tr}[\texttt{id}]\cdot \frac{2\pi i}{(\frac{n}{2}+1)!}\mathcal{A}_{11}dx'\nonumber\\
&+h'(0)\sum_{i=1}^{n-1}(a_{n}^{i})^{2}\cdot Vol(S^{n-2}){\rm tr}[\texttt{id}]\cdot \frac{2\pi i}{(\frac{n}{2}+1)!}\mathcal{A}_{12}dx'\nonumber\end{align}
\begin{align}
&+h'(0)\sum_{h=1}^{n}\sum_{i=1}^{n-1}(a_{h}^{i})^{2}\cdot \frac{1}{n-1}Vol(S^{n-2}){\rm tr}[\texttt{id}]\cdot \frac{2\pi i}{(\frac{n}{2}+2)!}\mathcal{A}_{13}dx'\nonumber\\
&+h'(0)\sum_{h=1}^{n}(a_{h}^{n})^{2}\cdot Vol(S^{n-2}){\rm tr}[\texttt{id}]\cdot \frac{2\pi i}{(\frac{n}{2}+2)!}\mathcal{A}_{14}dx',\nonumber
\end{align}
where
\begin{align}
\mathcal{A}_{11}
&=\Big[\frac{in(n-2)\xi_n^{2}-(n-2)\xi_n-i(n-2)}{8\left(\xi _n+i\right)^{\frac{n}{2}+1}}\Big]^{(\frac{n}{2}+1)}\bigg|_{\xi_n=i};\\
\mathcal{A}_{12}
&=\Big[\frac{-(n^{2}-4n+3)\xi_n^{3}-i(n-3)\xi_n^{2}+(3n-5)\xi_n+i}{8\left(\xi _n+i\right)^{\frac{n}{2}+1}}\Big]^{(\frac{n}{2}+1)}\bigg|_{\xi_n=i};\\
\mathcal{A}_{13}
&=\Big[\frac{-i(n^{2}-3n+2)\xi_n^{3}-2(n^{2}-3n+2)\xi_n^{2}-i(n^{2}-3n+2)\xi_n+2(n-2)}{8\left(\xi _n+i\right)^{\frac{n}{2}+1}}\Big]^{(\frac{n}{2}+2)}\bigg|_{\xi_n=i};\\
\mathcal{A}_{14}
&=\Big[\frac{-i(n^{2}-5n+6)\xi_n^{2}+2(n-2)\xi_n+i(n-2)}{8\left(\xi _n+i\right)^{\frac{n}{2}}}\Big]^{(\frac{n}{2}+2)}\bigg|_{\xi_n=i}.
\end{align}

\noindent {\bf case (b)}~$r=-1,l=-n+2,|\alpha|=j=k=0$.\\

It is easily seen that
\begin{align}
\Phi_4&=-i\int_{|\xi'|=1}\int^{+\infty}_{-\infty}{\rm tr}[\pi^{+}_{\xi_{n}}\sigma_{-1}({D}_{J}^{-1})
      \times\partial_{\xi_{n}}\sigma_{-n+2}({D}_{J}^{-n+3})](x_0)d\xi_n\sigma(\xi')dx'\\
&=i\int_{|\xi'|=1}\int^{+\infty}_{-\infty}{\rm tr}[\partial_{\xi_n}\pi^+_{\xi_n}\sigma_{-1}({D}_{J}^{-1})\times
\sigma_{-n+2}({D}_{J}^{-n+3})](x_0)d\xi_n\sigma(\xi')dx'.\nonumber
\end{align}

In the normal coordinate, $g^{ij}(x_{0})=\delta^{j}_{i}$ and $\partial_{x_{j}}(g^{\alpha\beta})(x_{0})=0$, if $j<n$; $\partial_{x_{j}}(g^{\alpha\beta})(x_{0})=h'(0)\delta^{\alpha}_{\beta}$, if $j=n$.
So by \cite{Wa3}, when $k<n$, we have $\Gamma^{n}(x_{0})=\frac{n-1}{2}h'(0)$ and $\Gamma^{k}(x_{0})=0.$
We thus get
\begin{align}
&\sigma_{-n+2}({D}_{J}^{-n+3})(x_{0})|_{|\xi'|=1}\\
&=-\frac{1}{4(1+\xi_{n}^{2})^{\frac{n}{2}-1}}h'(0)\sum_{\mu=1}^{n}\sum_{\nu=1}^{n-1}a_{\nu}^{\mu}c(dx_{\mu})c(dx_{n})c(dx_{\nu})\nonumber\\
&-\frac{n-2}{4(1+\xi_{n}^{2})^{\frac{n}{2}}}h'(0)\sum_{\omega=1}^{n}\sum_{\lambda,k=1}^{n-1}\xi_{k}\xi_{\lambda}a_{\omega}^{\lambda}c(dx_{k})c(dx_{n})c(dx_{\omega})\nonumber
\end{align}
\begin{align}
&-\frac{(n-2)\xi_{n}}{4(1+\xi_{n}^{2})^{\frac{n}{2}}}h'(0)\sum_{\omega=1}^{n}\sum_{k=1}^{n-1}\xi_{k}a_{\omega}^{n}c(dx_{k})c(dx_{n})c(dx_{\omega})\nonumber\\
&+\frac{(n^{2}-3n+2)\xi_{n}^{3}+(2n^{2}-5n+2)\xi_{n}}{4(1+\xi_{n}^{2})^{\frac{n}{2}+1}}h'(0)\sum_{\omega=1}^{n}\sum_{\lambda=1}^{n-1}\xi_{\lambda}a_{\omega}^{\lambda}c(dx_{\omega})\nonumber\\
&+\frac{(n^{2}-3n+2)\xi_{n}^{4}+(2n^{2}-5n+2)\xi_{n}^{2}}{4(1+\xi_{n}^{2})^{\frac{n}{2}+1}}h'(0)\sum_{\omega=1}^{n}a_{\omega}^{n}c(dx_{\omega})\nonumber\\
&+\frac{n-2}{2(1+\xi_{n}^{2})^{\frac{n}{2}}}\sum_{\alpha,\beta,\omega=1}^{n}\sum_{\lambda=1}^{n-1}\xi_{\lambda}a_{\alpha}^{\beta}a_{\omega}^{\lambda}c(dx_{\beta})c[(\nabla^{L}_{e_{\alpha}}J)(\xi^{*})]c(dx_{\omega})\nonumber\\
&+\frac{(n-2)\xi_{n}}{2(1+\xi_{n}^{2})^{\frac{n}{2}}}\sum_{\alpha,\beta,\omega=1}^{n}a_{\alpha}^{\beta}a_{\omega}^{n}c(dx_{\beta})c[(\nabla^{L}_{e_{\alpha}}J)(\xi^{*})]c(dx_{\omega})\nonumber\\
&-\frac{n-2}{(1+\xi_{n}^{2})^{\frac{n}{2}}}\sum_{h=1}^{n}\sum_{j,p=1}^{n-1}\xi_{j}\xi_{p}\partial_{x_{j}}(a_{h}^{p})c(dx_{h})\nonumber\\
&-\frac{(n-2)\xi_{n}}{(1+\xi_{n}^{2})^{\frac{n}{2}}}\sum_{h=1}^{n}\sum_{p=1}^{n-1}\xi_{p}\partial_{x_{n}}(a_{h}^{p})c(dx_{h})\nonumber\\
&-\frac{(n-2)\xi_{n}}{(1+\xi_{n}^{2})^{\frac{n}{2}}}\sum_{h=1}^{n}\sum_{j=1}^{n-1}\xi_{j}\partial_{x_{j}}(a_{h}^{n})c(dx_{h})\nonumber\\
&-\frac{(n-2)\xi_{n}^{2}}{(1+\xi_{n}^{2})^{\frac{n}{2}}}\sum_{h=1}^{n}\partial_{x_{n}}(a_{h}^{n})c(dx_{h})\nonumber\\
&-\frac{(n-2)\xi_{n}}{(1+\xi_{n}^{2})^{\frac{n}{2}}}\sum_{p,h=1}^{n-1}\xi_{p}a_{h}^{p}\partial_{x_{n}}(c(dx_{h}))\nonumber\\
&-\frac{(n-2)\xi_{n}^{2}}{(1+\xi_{n}^{2})^{\frac{n}{2}}}\sum_{h=1}^{n-1}a_{h}^{n}\partial_{x_{n}}(c(dx_{h})).\nonumber
\end{align}

Summarizing, we have
\begin{align}
\Phi_4
&=h'(0)\sum_{i=1}^{n-1}(a_{n}^{i})^{2}\cdot \frac{1}{n-1}Vol(S^{n-2}){\rm tr}[\texttt{id}]\cdot \frac{2\pi i}{(\frac{n}{2}+1)!}\mathcal{B}_{0}dx'\\
&+h'(0)\sum_{i=1}^{n-1}(a_{n}^{i})^{2}\cdot Vol(S^{n-2}){\rm tr}[\texttt{id}]\cdot \frac{2\pi i}{(\frac{n}{2})!}\mathcal{B}_{1}dx'\nonumber\\
&+h'(0)\sum_{i=1}^{n-1}a_{i}^{i}a_{n}^{n}\cdot \frac{1}{n-1}Vol(S^{n-2}){\rm tr}[\texttt{id}]\cdot \frac{2\pi i}{(\frac{n}{2}+1)!}\cdot-\mathcal{B}_{0}dx'\nonumber
\end{align}
\begin{align}
&+h'(0)\sum_{i=1}^{n-1}a_{i}^{i}a_{n}^{n}\cdot Vol(S^{n-2}){\rm tr}[\texttt{id}]\cdot \frac{2\pi i}{(\frac{n}{2})!}\cdot-\mathcal{B}_{1}dx'\nonumber\\
&+h'(0)\sum_{l=1}^{n}(a_{n}^{l})^{2}\cdot Vol(S^{n-2}){\rm tr}[\texttt{id}]\cdot \frac{2\pi i}{(\frac{n}{2}+2)!}\mathcal{B}_{2}dx'\nonumber\\
&+h'(0)\sum_{i=1}^{n-1}(a_{n}^{i})^{2}\cdot Vol(S^{n-2}){\rm tr}[\texttt{id}]\cdot \frac{2\pi i}{(\frac{n}{2}+1)!}\mathcal{B}_{3}dx'\nonumber\\
&+h'(0)\sum_{l=1}^{n}\sum_{i=1}^{n-1}(a_{l}^{i})^{2}\cdot \frac{1}{n-1}Vol(S^{n-2}){\rm tr}[\texttt{id}]\cdot \frac{2\pi i}{(\frac{n}{2}+2)!}\mathcal{B}_{4}dx'\nonumber\\
&+h'(0)\sum_{i,j=1}^{n-1}(a_{j}^{i})^{2}\cdot \frac{1}{n-1}Vol(S^{n-2}){\rm tr}[\texttt{id}]\cdot \frac{2\pi i}{(\frac{n}{2}+1)!}\mathcal{B}_{5}dx'\nonumber\\
&+\sum_{l,\alpha=1}^{n}\sum_{i=1}^{n-1}(a_{l}^{i})^{2}g^{M}(J(dx_{\alpha}), (\nabla^{L}_{e_{\alpha}}J)e_{n})\cdot \frac{1}{n-1}Vol(S^{n-2}){\rm tr}[\texttt{id}]\cdot \frac{2\pi i}{(\frac{n}{2}+1)!}\mathcal{B}_{5}dx'\nonumber\\
&+\sum_{i=1}^{n-1}g^{M}(J(dx_{i}), (\nabla^{L}_{e_{n}}J)e_{i})\cdot \frac{1}{n-1}Vol(S^{n-2}){\rm tr}[\texttt{id}]\cdot \frac{2\pi i}{(\frac{n}{2}+1)!}\cdot-2\mathcal{B}_{0}dx'\nonumber\\
&+\sum_{i=1}^{n-1}g^{M}(J(dx_{n}), (\nabla^{L}_{e_{i}}J)e_{i})\cdot \frac{1}{n-1}Vol(S^{n-2}){\rm tr}[\texttt{id}]\cdot \frac{2\pi i}{(\frac{n}{2}+1)!}\cdot2\mathcal{B}_{0}dx'\nonumber\\
&+\sum_{l,\alpha=1}^{n}(a_{l}^{n})^{2}g^{M}(J(dx_{\alpha}), (\nabla^{L}_{e_{\alpha}}J)e_{n})\cdot Vol(S^{n-2}){\rm tr}[\texttt{id}]\cdot \frac{2\pi i}{(\frac{n}{2}+1)!}\mathcal{B}_{3}dx'\nonumber\\
&+\sum_{l=1}^{n}\sum_{i=1}^{n-1}a_{l}^{i}\partial_{x_{i}}(a_{l}^{n})\cdot \frac{1}{n-1}Vol(S^{n-2}){\rm tr}[\texttt{id}]\cdot \frac{2\pi i}{(\frac{n}{2}+1)!}\cdot4\mathcal{B}_{0}dx',\nonumber
\end{align}
where
\begin{align}
\mathcal{B}_{0}
&=\Big[-\frac{i(n-2)\xi_n+n-2}{8\left(\xi _n+i\right)^{\frac{n}{2}}}\Big]^{(\frac{n}{2}+1)}\bigg|_{\xi_n=i};\\
\mathcal{B}_{1}
&=\Big[-\frac{1}{8\left(\xi _n+i\right)^{\frac{n}{2}-1}}\Big]^{(\frac{n}{2})}\bigg|_{\xi_n=i};\\
\mathcal{B}_{2}
&=\Big[-\frac{(n^{2}-3n+2)\xi_n^{4}+(2n^{2}-5n+2)\xi_n^{2}}{8\left(\xi _n+i\right)^{\frac{n}{2}+1}}\Big]^{(\frac{n}{2}+2)}\bigg|_{\xi_n=i};\\
\mathcal{B}_{3}
&=\Big[\frac{(n-2)\xi_n^{2}}{4\left(\xi _n+i\right)^{\frac{n}{2}}}\Big]^{(\frac{n}{2}+1)}\bigg|_{\xi_n=i};
\end{align}
\begin{align}
\mathcal{B}_{4}
&=\Big[\frac{i(n^{2}-3n+2)\xi_n^{3}+i(2n^{2}-5n+2)\xi_n}{8\left(\xi _n+i\right)^{\frac{n}{2}+1}}\Big]^{(\frac{n}{2}+2)}\bigg|_{\xi_n=i};\\
\mathcal{B}_{5}
&=\Big[-\frac{i(n-2)\xi_n}{4\left(\xi _n+i\right)^{\frac{n}{2}}}\Big]^{(\frac{n}{2}+1)}\bigg|_{\xi_n=i}.
\end{align}

\noindent {\bf  case (c)}~$r=-2,l=-n+3,|\alpha|=j=k=0$.\\

We calculate
\begin{equation}
\Phi_5=-i\int_{|\xi'|=1}\int^{+\infty}_{-\infty}{\rm tr}[\pi^{+}_{\xi_{n}}\sigma_{-2}({D}_{J}^{-1})
      \times\partial_{\xi_{n}}\sigma_{-n+3}({D}_{J}^{-n+3})](x_0)d\xi_n\sigma(\xi')dx'.
\end{equation}

We follow the notation of \cite{LSW1}.
\begin{align}
A_1(x_0)&=\frac{c[J(\xi)]\sigma_{0}({D}_{J})(x_0)c[J(\xi)]}{(1+\xi_n^2)^2};\\
A_2(x_0)&=\frac{c[J(\xi)]}{(1+\xi_n^2)^2}\Big[\sum_{j,p,h=1}^{n}\xi_p\partial_{x_j}(a_{h}^{p})c[J(dx_j)]c(dx_h)+\sum_{p=1}^{n}\sum_{h=1}^{n-1}\xi_pa_{h}^{p}c[J(dx_n)]\partial_{x_n}(c(dx_h))\Big];\\
A_3(x_0)&=\frac{c[J(\xi)]}{(1+\xi_n^2)^3}c[J(dx_n)]c[J(\xi)],
\end{align}
means that
\begin{align}
&\pi^+_{\xi_n}\sigma_{-2}({D}_{J}^{-1})(x_0)|_{|\xi'|=1}=\pi^+_{\xi_n}(A_1(x_0))+\pi^+_{\xi_n}(A_2(x_0))-h'(0)\pi^+_{\xi_n}(A_3(x_0)).
\end{align}

Computations show that
\begin{align}
\pi^+_{\xi_n}(A_1(x_0))
&=\frac{i\xi_n}{16(\xi_n-i)^2}h'(0)\sum_{l,\gamma,\mu=1}^{n}\sum_{\nu=1}^{n-1}a_{l}^{n}a_{\gamma}^{n}a_{\nu}^{\mu}c(dx_l)c(dx_\mu)c(dx_n)c(dx_\nu)c(dx_{\gamma})\\
&+\frac{i}{16(\xi_n-i)^2}h'(0)\sum_{l,\gamma,\mu=1}^{n}\sum_{q,\nu=1}^{n-1}\xi_{q}a_{l}^{q}a_{\gamma}^{n}a_{\nu}^{\mu}c(dx_l)c(dx_\mu)c(dx_n)c(dx_\nu)c(dx_{\gamma})\nonumber\\
&+\frac{i}{16(\xi_n-i)^2}h'(0)\sum_{l,\gamma,\mu=1}^{n}\sum_{\alpha,\nu=1}^{n-1}\xi_{\alpha}a_{l}^{n}a_{\gamma}^{\alpha}a_{\nu}^{\mu}c(dx_l)c(dx_\mu)c(dx_n)c(dx_\nu)c(dx_{\gamma})\nonumber\\
&+\frac{i\xi_n+2}{16(\xi_n-i)^2}h'(0)\sum_{l,\gamma,\mu=1}^{n}\sum_{q,\alpha,\nu=1}^{n-1}\xi_{q}\xi_{\alpha}a_{l}^{q}a_{\gamma}^{\alpha}a_{\nu}^{\mu}c(dx_l)c(dx_\mu)c(dx_n)c(dx_\nu)c(dx_{\gamma}).\nonumber
\end{align}
By computation, we have
\begin{align}
\partial_{\xi_n}\left(\frac{ic[J(\xi)]}{|\xi|^{n-2}}\right)(x_0)|_{|\xi'|=1}
&=i\sum^{n}_{\beta=1}\sum^{n-1}_{i=1}\xi_{i}a_{\beta}^ic(dx_{\beta})\partial_{\xi_n}\left(\frac{1}{(1+\xi_{n}^{2})^{\frac{n}{2}-1}}\right)\\
&+i\sum^{n}_{\beta=1}a_{\beta}^nc(dx_{\beta})\partial_{\xi_n}\left(\frac{\xi_{n}}{(1+\xi_{n}^{2})^{\frac{n}{2}-1}}\right)\nonumber\\
&=-\frac{i(n-2)\xi_n}{\left(1+\xi _n^2\right)^{\frac{n}{2}}}\sum^{n}_{\beta=1}\sum^{n-1}_{i=1}\xi_{i}a_{\beta}^ic(dx_{\beta})\nonumber\\
&+\frac{i[1-(n-3)\xi_n^2]}{\left(1+\xi _n^2\right)^{\frac{n}{2}}}\sum^{n}_{\beta=1}a_{\beta}^nc(dx_{\beta}).\nonumber
\end{align}

By $\int_{|\xi'|=1}{\{\xi_{i_1}\cdot\cdot\cdot\xi_{i_{2d+1}}}\}\sigma(\xi')=0$ and $\int_{|\xi'|=1}\sum_{i,j=1}^{n-1}\xi_i\xi_j\sigma(\xi')=\frac{1}{n-1}\sum_{i,j=1}^{n-1}\delta_{j}^{i}Vol(S^{n-2})$, then we have
\begin{align}
&-i\int_{|\xi'|=1}\int^{+\infty}_{-\infty}{\rm tr} [\pi^+_{\xi_n}(A_1(x_0)) \times \partial_{\xi_n}\sigma_{-n+3}({D}_{J}^{-n+3})](x_0)d\xi_n\sigma(\xi')dx'\\
&=-i\int^{+\infty}_{-\infty}\frac{(n-3)\xi_n^{3}-\xi_n}{16 \left(\xi_n-i\right)^{\frac{n}{2}+2} \left(\xi _n+i\right)^{\frac{n}{2}}}h'(0)\nonumber\\
&\cdot \sum_{l=1}^{n}\sum_{\nu=1}^{n-1}(-(a_{\nu}^{n})^2(a_{l}^{n})^2+(a_{l}^{n})^2a_{\nu}^{\nu}a_{n}^{n})\cdot Vol(S^{n-2}){\rm tr}[\texttt{id}]d\xi_ndx'\nonumber\\
&-i\int^{+\infty}_{-\infty}\frac{(n-2)\xi_n}{16 \left(\xi_n-i\right)^{\frac{n}{2}+2} \left(\xi _n+i\right)^{\frac{n}{2}}}h'(0)\nonumber\\
&\cdot \sum_{l=1}^{n}\sum_{\nu,i=1}^{n-1}(-(a_{\nu}^{n})^2(a_{l}^{i})^2+(a_{l}^{i})^2a_{\nu}^{\nu}a_{n}^{n})\cdot \frac{2}{n-1}Vol(S^{n-2}){\rm tr}[\texttt{id}]d\xi_ndx'\nonumber\\
&-i\int^{+\infty}_{-\infty}\frac{(n-3)\xi_n^{3}-2i(n-3)\xi_n^{2}-\xi_n+2i}{16 \left(\xi_n-i\right)^{\frac{n}{2}+2} \left(\xi _n+i\right)^{\frac{n}{2}}}h'(0)\nonumber\\
&\cdot \sum_{l=1}^{n}\sum_{\nu,i=1}^{n-1}((a_{\nu}^{n})^2(a_{l}^{i})^2-(a_{l}^{i})^2a_{\nu}^{\nu}a_{n}^{n})\cdot \frac{1}{n-1}Vol(S^{n-2}){\rm tr}[\texttt{id}]d\xi_ndx'.\nonumber
\end{align}

\cite{LSW1} also shown that
\begin{align}
\pi^+_{\xi_n}(A_2(x_0))
&=-\frac{i\xi_n}{4(\xi_n-i)^2}\sum_{l,j,h,y=1}^{n}a_{l}^{n}a_{y}^{j}\partial_{x_j}(a_{h}^{n})c(dx_{l})c(dx_{y})c(dx_{h})\\
&-\frac{i}{4(\xi_n-i)^2}\sum_{l,j,h,y=1}^{n}\sum_{q=1}^{n-1}\xi_{q}a_{l}^{q}a_{y}^{j}\partial_{x_j}(a_{h}^{n})c(dx_{l})c(dx_{y})c(dx_{h})\nonumber
\end{align}
\begin{align}
&-\frac{i}{4(\xi_n-i)^2}\sum_{l,j,h,y=1}^{n}\sum_{p=1}^{n-1}\xi_{p}a_{l}^{n}a_{y}^{j}\partial_{x_j}(a_{h}^{p})c(dx_{l})c(dx_{y})c(dx_{h})\nonumber\\
&-\frac{i\xi_n+2}{4(\xi_n-i)^2}\sum_{l,j,h,y=1}^{n}\sum_{q,p=1}^{n-1}\xi_{q}\xi_{p}a_{l}^{q}a_{y}^{j}\partial_{x_j}(a_{h}^{p})c(dx_{l})c(dx_{y})c(dx_{h})\nonumber\\
&-\frac{i\xi_n}{8(\xi_n-i)^2}h'(0)\sum_{l,z=1}^{n}\sum_{h=1}^{n-1}a_{l}^{n}a_{h}^{n}a_{z}^{n}c(dx_{l})c(dx_{z})c(dx_{h})\nonumber\\
&-\frac{i}{8(\xi_n-i)^2}h'(0)\sum_{l,z=1}^{n}\sum_{q,h=1}^{n-1}\xi_{q}a_{l}^{q}a_{h}^{n}a_{z}^{n}c(dx_{l})c(dx_{z})c(dx_{h})\nonumber\\
&-\frac{i}{8(\xi_n-i)^2}h'(0)\sum_{l,z=1}^{n}\sum_{p,h=1}^{n-1}\xi_{p}a_{l}^{n}a_{h}^{p}a_{z}^{n}c(dx_{l})c(dx_{z})c(dx_{h})\nonumber\\
&-\frac{i\xi_n+2}{8(\xi_n-i)^2}h'(0)\sum_{l,z=1}^{n}\sum_{q,p,h=1}^{n-1}\xi_{q}\xi_{p}a_{l}^{q}a_{h}^{p}a_{z}^{n}c(dx_{l})c(dx_{z})c(dx_{h}).\nonumber
\end{align}
Accordingly, we have
\begin{align}
&-i\int_{|\xi'|=1}\int^{+\infty}_{-\infty}{\rm tr} [\pi^+_{\xi_n}(A_2(x_0)) \times \partial_{\xi_n}\sigma_{-n+3}({D}_{J}^{-n+3})](x_0)d\xi_n\sigma(\xi')dx'\\
&=-i\int^{+\infty}_{-\infty}\frac{-(n-3)\xi_n^{3}+\xi_n}{4 \left(\xi_n-i\right)^{\frac{n}{2}+2} \left(\xi _n+i\right)^{\frac{n}{2}}}\sum_{j,\beta,\eta=1}^{n}(a_{\beta}^{n})^2a_{\eta}^{j}\partial_{x_j}(a_{\eta}^{n})\cdot Vol(S^{n-2}){\rm tr}[\texttt{id}]d\xi_ndx'\nonumber\\
&-i\int^{+\infty}_{-\infty}-\frac{(n-2)\xi_n}{4 \left(\xi_n-i\right)^{\frac{n}{2}+2} \left(\xi _n+i\right)^{\frac{n}{2}}}\sum_{j,\beta,\eta=1}^{n}\sum_{i=1}^{n-1}(a_{\beta}^{i})^2a_{\eta}^{j}\partial_{x_j}(a_{\eta}^{n})\cdot\frac{1}{n-1} Vol(S^{n-2}){\rm tr}[\texttt{id}]d\xi_ndx'\nonumber\\
&-i\int^{+\infty}_{-\infty}\frac{(n-2)\xi_n}{4 \left(\xi_n-i\right)^{\frac{n}{2}+2} \left(\xi _n+i\right)^{\frac{n}{2}}}\sum_{l=1}^{n}\sum_{i=1}^{n-1}a_{l}^{n}\partial_{x_i}(a_{l}^{i})\cdot\frac{1}{n-1} Vol(S^{n-2}){\rm tr}[\texttt{id}]d\xi_ndx'\nonumber\\
&-i\int^{+\infty}_{-\infty}-\frac{(n-3)\xi_n^{3}-2i(n-3)\xi_n^{2}-\xi_n+2i}{4 \left(\xi_n-i\right)^{\frac{n}{2}+2} \left(\xi _n+i\right)^{\frac{n}{2}}}\sum_{l=1}^{n}\sum_{i=1}^{n-1}a_{l}^{n}\partial_{x_i}(a_{l}^{i})\cdot\frac{1}{n-1} Vol(S^{n-2}){\rm tr}[\texttt{id}]d\xi_ndx'\nonumber\\
&-i\int^{+\infty}_{-\infty}\frac{-(n-3)\xi_n^{3}+\xi_n}{8 \left(\xi_n-i\right)^{\frac{n}{2}+2} \left(\xi _n+i\right)^{\frac{n}{2}}}h'(0)\sum_{l=1}^{n}\sum_{\nu=1}^{n-1}(a_{\nu}^{n})^2(a_{l}^{n})^2\cdot Vol(S^{n-2}){\rm tr}[\texttt{id}]d\xi_ndx'\nonumber\\
&-i\int^{+\infty}_{-\infty}-\frac{(n-2)\xi_n}{8 \left(\xi_n-i\right)^{\frac{n}{2}+2} \left(\xi _n+i\right)^{\frac{n}{2}}}h'(0)\sum_{l=1}^{n}\sum_{\nu,i=1}^{n-1}(a_{\nu}^{n})^2(a_{l}^{i})^2\cdot\frac{1}{n-1} Vol(S^{n-2}){\rm tr}[\texttt{id}]d\xi_ndx'\nonumber\end{align}
\begin{align}
&-i\int^{+\infty}_{-\infty}-\frac{(n-2)\xi_n}{8 \left(\xi_n-i\right)^{\frac{n}{2}+2} \left(\xi _n+i\right)^{\frac{n}{2}}}h'(0)\sum_{l=1}^{n}\sum_{\nu,i=1}^{n-1}(a_{\nu}^{i})^2(a_{l}^{n})^2\cdot\frac{1}{n-1} Vol(S^{n-2}){\rm tr}[\texttt{id}]d\xi_ndx'\nonumber\\
&-i\int^{+\infty}_{-\infty}\frac{(n-3)\xi_n^{3}-2i(n-3)\xi_n^{2}-\xi_n+2i}{8 \left(\xi_n-i\right)^{\frac{n}{2}+2} \left(\xi _n+i\right)^{\frac{n}{2}}}h'(0)\sum_{l=1}^{n}\sum_{\nu,i=1}^{n-1}(a_{\nu}^{i})^2(a_{l}^{n})^2\cdot\frac{1}{n-1} Vol(S^{n-2}){\rm tr}[\texttt{id}]d\xi_ndx'.\nonumber
\end{align}

Since
\begin{align}
-h'(0)\pi^+_{\xi_n}(A_3(x_0))
&=\frac{3\xi_n+i\xi_n^2}{16(\xi_n-i)^3}h'(0)\sum_{l,w,\gamma=1}^{n}a_{l}^{n}a_{w}^{n}a_{\gamma}^{n}c(dx_{l})c(dx_{w})c(dx_{\gamma})\\
&+\frac{i\xi_n+3}{16(\xi_n-i)^3}h'(0)\sum_{l,w,\gamma=1}^{n}\sum_{q=1}^{n-1}\xi_{q}a_{l}^{q}a_{w}^{n}a_{\gamma}^{n}c(dx_{l})c(dx_{w})c(dx_{\gamma})\nonumber\\
&+\frac{i\xi_n+3}{16(\xi_n-i)^3}h'(0)\sum_{l,w,\gamma=1}^{n}\sum_{\alpha=1}^{n-1}\xi_{\alpha}a_{l}^{n}a_{w}^{n}a_{\gamma}^{\alpha}c(dx_{l})c(dx_{w})c(dx_{\gamma})\nonumber\\
&+\frac{-8i+9\xi_n+3i\xi_n^2}{16(\xi_n-i)^3}h'(0)\sum_{l,w,\gamma=1}^{n}\sum_{q,\alpha=1}^{n-1}\xi_{q}\xi_{\alpha}a_{l}^{q}a_{w}^{n}a_{\gamma}^{\alpha}c(dx_{l})c(dx_{w})c(dx_{\gamma}),\nonumber
\end{align}
it is sufficient to show that
\begin{align}
&-i\int_{|\xi'|=1}\int^{+\infty}_{-\infty}{\rm tr} [-h'(0)\pi^+_{\xi_n}(A_3(x_0)) \times \partial_{\xi_n}\sigma_{-n+3}({D}_{J}^{-n+3})](x_0)d\xi_n\sigma(\xi')dx'\\
&=-i\int^{+\infty}_{-\infty}\frac{(n-3)\xi_n^{4}-3i(n-3)\xi_n^{3}-\xi_n^{2}+3i\xi_n}{16 \left(\xi_n-i\right)^{\frac{n}{2}+3} \left(\xi _n+i\right)^{\frac{n}{2}}}h'(0)\sum_{\beta,l=1}^{n}(a_{\beta}^{n})^2(a_{l}^{n})^2\cdot Vol(S^{n-2}){\rm tr}[\texttt{id}]d\xi_ndx'\nonumber\\
&-i\int^{+\infty}_{-\infty}\frac{(n-2)\xi_n^{2}-3i(n-2)\xi_n}{16 \left(\xi_n-i\right)^{\frac{n}{2}+3} \left(\xi _n+i\right)^{\frac{n}{2}}}h'(0)\sum_{\beta,\omega=1}^{n}\sum_{i=1}^{n-1}(a_{\beta}^{i})^2(a_{\omega}^{n})^2\cdot \frac{2}{n-1} Vol(S^{n-2}){\rm tr}[\texttt{id}]d\xi_ndx'\nonumber\\
&-i\int^{+\infty}_{-\infty}-\frac{3(n-3)\xi_n^{4}-9i(n-3)\xi_n^{3}-(8n-21)\xi_n^{2}+9i\xi_n+8}{16 \left(\xi_n-i\right)^{\frac{n}{2}+3} \left(\xi _n+i\right)^{\frac{n}{2}}}h'(0)\sum_{\beta,l=1}^{n}\sum_{i=1}^{n-1}(a_{l}^{i})^2(a_{\beta}^{n})^2\nonumber\\
&\cdot \frac{1}{n-1} Vol(S^{n-2}){\rm tr}[\texttt{id}]d\xi_ndx'.\nonumber
\end{align}

It follows immediately that
\begin{align}
\sum_{l=1}^{n}\sum_{i=1}^{n-1}a_{l}^{n}\partial_{x_i}(a_{l}^{i})&=\sum_{l,i=1}^{n}a_{l}^{n}\partial_{x_i}(a_{l}^{i})-\sum_{l=1}^{n}a_{l}^{n}\partial_{x_n}(a_{l}^{n})
\end{align}
\begin{align}
&=\sum_{l,i=1}^{n}a_{l}^{n}\partial_{x_i}(a_{l}^{i})-\partial_{x_n}(\delta_{n}^{n})=\sum_{l,i=1}^{n}a_{l}^{n}\partial_{x_i}(a_{l}^{i}).\nonumber
\end{align}

On account of the above result,
\begin{align}
\Phi_5&=\sum_{j,\beta,\eta=1}^{n}(a_{\beta}^{n})^2a_{\eta}^{j}\partial_{x_j}(a_{\eta}^{n})\cdot Vol(S^{n-2}){\rm tr}[\texttt{id}]\cdot \frac{2\pi i}{(\frac{n}{2}+1)!}\mathcal{C}_{0}dx'\\
&+\sum_{j,\beta,\eta=1}^{n}\sum_{i=1}^{n-1}(a_{\beta}^{i})^2a_{\eta}^{j}\partial_{x_j}(a_{\eta}^{n})\cdot \frac{1}{n-1}Vol(S^{n-2}){\rm tr}[\texttt{id}]\cdot \frac{2\pi i}{(\frac{n}{2}+1)!}\mathcal{C}_{1}dx'\nonumber\\
&+h'(0)\sum_{l=1}^{n}\sum_{\nu=1}^{n-1}(a_{l}^{n})^2a_{\nu}^{\nu}a_{n}^{n}\cdot Vol(S^{n-2}){\rm tr}[\texttt{id}]\cdot \frac{2\pi i}{(\frac{n}{2}+1)!}\cdot-\frac{1}{4}\mathcal{C}_{0}dx'\nonumber\\
&+h'(0)\sum_{l=1}^{n}\sum_{\nu,i=1}^{n-1}(a_{l}^{i})^2a_{\nu}^{\nu}a_{n}^{n}\cdot \frac{1}{n-1}Vol(S^{n-2}){\rm tr}[\texttt{id}]\cdot \frac{2\pi i}{(\frac{n}{2}+1)!}\mathcal{C}_{2}dx'\nonumber\\
&+\sum_{l,j=1}^{n}a_{l}^{j}\partial_{x_j}(a_{l}^{n})\cdot \frac{1}{n-1}Vol(S^{n-2}){\rm tr}[\texttt{id}]\cdot \frac{2\pi i}{(\frac{n}{2}+1)!}\mathcal{C}_{3}dx'\nonumber\\
&+h'(0)\sum_{\beta,l=1}^{n}(a_{\beta}^{n})^2(a_{l}^{n})^2\cdot Vol(S^{n-2}){\rm tr}[\texttt{id}]\cdot \frac{2\pi i}{(\frac{n}{2}+2)!}\mathcal{C}_{4}dx'\nonumber\\
&+h'(0)\sum_{\beta,l=1}^{n}\sum_{i=1}^{n-1}(a_{\beta}^{n})^2(a_{l}^{i})^2\cdot \frac{1}{n-1}Vol(S^{n-2}){\rm tr}[\texttt{id}]\cdot \frac{2\pi i}{(\frac{n}{2}+2)!}\mathcal{C}_{5}dx'\nonumber\\
&+h'(0)\sum_{l=1}^{n}\sum_{\nu=1}^{n-1}(a_{\nu}^{n})^2(a_{l}^{n})^2\cdot Vol(S^{n-2}){\rm tr}[\texttt{id}]\cdot \frac{2\pi i}{(\frac{n}{2}+1)!}\cdot\frac{3}{4}\mathcal{C}_{0}dx'\nonumber\\
&+h'(0)\sum_{l=1}^{n}\sum_{\nu,i=1}^{n-1}(a_{\nu}^{n})^2(a_{l}^{i})^2\cdot \frac{1}{n-1}Vol(S^{n-2}){\rm tr}[\texttt{id}]\cdot \frac{2\pi i}{(\frac{n}{2}+1)!}\mathcal{C}_{6}dx'\nonumber\\
&+h'(0)\sum_{l=1}^{n}\sum_{\nu,i=1}^{n-1}(a_{\nu}^{i})^2(a_{l}^{n})^2\cdot \frac{1}{n-1}Vol(S^{n-2}){\rm tr}[\texttt{id}]\cdot \frac{2\pi i}{(\frac{n}{2}+1)!}\cdot\frac{1}{2}\mathcal{C}_{3}dx',\nonumber
\end{align}
where
\begin{align}
\mathcal{C}_{0}
&=\Big[\frac{i(n-3)\xi_n^{3}-i\xi_n}{4\left(\xi _n+i\right)^{\frac{n}{2}}}\Big]^{(\frac{n}{2}+1)}\bigg|_{\xi_n=i};\\
\mathcal{C}_{1}
&=\Big[\frac{i(n-2)\xi_n}{4\left(\xi _n+i\right)^{\frac{n}{2}}}\Big]^{(\frac{n}{2}+1)}\bigg|_{\xi_n=i};\\
\mathcal{C}_{2}
&=\Big[\frac{i(n-3)\xi_n^{3}+2(n-3)\xi_n^{2}-i(2n-3)\xi_n-2}{16\left(\xi _n+i\right)^{\frac{n}{2}}}\Big]^{(\frac{n}{2}+1)}\bigg|_{\xi_n=i};
\end{align}
\begin{align}
\mathcal{C}_{3}
&=\Big[\frac{-i(n-3)\xi_n^{3}-2(n-3)\xi_n^{2}+i(n-1)\xi_n+2}{4\left(\xi _n+i\right)^{\frac{n}{2}}}\Big]^{(\frac{n}{2}+1)}\bigg|_{\xi_n=i};\\
\mathcal{C}_{4}
&=\Big[\frac{-i(n-3)\xi_n^{4}-3(n-3)\xi_n^{3}+i\xi_n^{2}+3\xi_n}{16\left(\xi _n+i\right)^{\frac{n}{2}}}\Big]^{(\frac{n}{2}+2)}\bigg|_{\xi_n=i};\\
\mathcal{C}_{5}
&=\Big[\frac{3i(n-3)\xi_n^{4}+9(n-3)\xi_n^{3}-5i(2n-5)\xi_n^{2}-3(2n-1)\xi_n+8i}{16\left(\xi _n+i\right)^{\frac{n}{2}}}\Big]^{(\frac{n}{2}+2)}\bigg|_{\xi_n=i};\\
\mathcal{C}_{6}
&=\Big[\frac{-i(n-3)\xi_n^{3}-2(n-3)\xi_n^{2}-i(-4n+7)\xi_n+2}{16\left(\xi _n+i\right)^{\frac{n}{2}}}\Big]^{(\frac{n}{2}+1)}\bigg|_{\xi_n=i}.
\end{align}

A simple calculation shows that
\begin{align}
\sum_{i=1}^{n-1}g^{M}(J(dx_{i}), (\nabla^{L}_{e_{n}}J)e_{i})(x_{0})&=\sum_{i=1}^{n-1}g^{M}(J(e_{i}), (\nabla^{L}_{e_{n}}J)e_{i})(x_{0})\\
&=\sum_{i=1}^{n-1}g^{M}(J(e_{i}), \nabla^{L}_{e_{n}}J(e_{i})-J(\nabla^{L}_{e_{n}}e_{i}))(x_{0})\nonumber\\
&=\sum_{i=1}^{n-1}g^{M}(J(e_{i}), \nabla^{L}_{e_{n}}J(e_{i}))(x_{0})-\sum_{i=1}^{n-1}g^{M}(J(e_{i}), J(\nabla^{L}_{e_{n}}e_{i}))(x_{0})\nonumber\\
&=\sum_{i=1}^{n-1}g^{M}(J(e_{i}), \nabla^{L}_{e_{n}}J(e_{i}))(x_{0})-\sum_{i=1}^{n-1}g^{M}(e_{i}, \nabla^{L}_{e_{n}}e_{i})(x_{0})\nonumber\\
&=\frac{1}{2}\sum_{i=1}^{n-1}e_{n}(g^{M}(J(e_{i}), J(e_{i})))(x_{0})-\frac{1}{2}\sum_{i=1}^{n-1}e_{n}(g^{M}(e_{i}, e_{i}))(x_{0})=0.\nonumber
\end{align}
Since
\begin{align}
&\sum_{i=1}^{n-1}g^{M}(J(dx_{n}), (\nabla^{L}_{e_{i}}J)e_{i})(x_{0})+\sum_{i=1}^{n-1}g^{M}(J(dx_{i}), (\nabla^{L}_{e_{i}}J)e_{n})(x_{0})\\
&=\sum_{i=1}^{n-1}g^{M}(J(e_{n}), \nabla^{L}_{e_{i}}J(e_{i}))(x_{0})-\sum_{i=1}^{n-1}g^{M}(J(e_{n}), J(\nabla^{L}_{e_{i}}e_{i}))(x_{0})\nonumber\\
&+\sum_{i=1}^{n-1}g^{M}(J(e_{i}), \nabla^{L}_{e_{i}}J(e_{n}))(x_{0})-\sum_{i=1}^{n-1}g^{M}(J(e_{i}), J(\nabla^{L}_{e_{i}}e_{n}))(x_{0})=0,\nonumber
\end{align}
we get 
\begin{align}
\sum_{i=1}^{n-1}g^{M}(J(dx_{n}), (\nabla^{L}_{e_{i}}J)e_{i})(x_{0})=-\sum_{i=1}^{n-1}g^{M}(J(e_{i}), (\nabla^{L}_{e_{i}}J)e_{n})(x_{0}).
\end{align}
Similarly, we have
\begin{align}
&\sum_{\alpha=1}^{n}g^{M}(J(dx_{\alpha}), (\nabla^{L}_{e_{\alpha}}J)e_{n})(x_{0})=\sum_{\alpha=1}^{n-1}g^{M}(J(e_{\alpha}), (\nabla^{L}_{e_{\alpha}}J)e_{n})(x_{0}).
\end{align}

In summary,
\begin{align}
\Phi
&=\Phi_1+\Phi_2+\Phi_3+\Phi_4+\Phi_5\\
&=h'(0)\cdot Vol(S^{n-2}){\rm tr}[\texttt{id}]\cdot \frac{2\pi i}{(\frac{n}{2}+2)!}\mathcal{D}_{0}dx'\nonumber\\
&+\sum_{h=1}^{n}\sum_{i=1}^{n-1}a_{h}^{i}\partial_{x_i}(a_{h}^{n})\cdot Vol(S^{n-2}){\rm tr}[\texttt{id}]\cdot \frac{2\pi i}{(\frac{n}{2})!}\mathcal{D}_{1}dx'\nonumber\\
&+\sum_{h=1}^{n}\sum_{i=1}^{n-1}a_{h}^{i}\partial_{x_i}(a_{h}^{n})\cdot Vol(S^{n-2}){\rm tr}[\texttt{id}]\cdot \frac{2\pi i}{(\frac{n}{2}+1)!}\mathcal{D}_{2}dx'\nonumber\\
&+h'(0)\sum_{\nu,i=1}^{n-1}(a_{\nu}^{i})^2\cdot Vol(S^{n-2}){\rm tr}[\texttt{id}]\cdot \frac{2\pi i}{(\frac{n}{2}+1)!}\mathcal{D}_{3}dx'\nonumber\\
&+h'(0)\sum_{i=1}^{n-1}(a_{n}^{i})^2\cdot Vol(S^{n-2}){\rm tr}[\texttt{id}]\cdot \frac{2\pi i}{(\frac{n}{2})!}\mathcal{D}_{4}dx'\nonumber\\
&+h'(0)\sum_{i=1}^{n-1}(a_{n}^{i})^2\cdot Vol(S^{n-2}){\rm tr}[\texttt{id}]\cdot \frac{2\pi i}{(\frac{n}{2}+1)!}\mathcal{D}_{5}dx'\nonumber\\
&+h'(0)\sum_{i=1}^{n-1}a_{i}^{i}a_{n}^{n}\cdot Vol(S^{n-2}){\rm tr}[\texttt{id}]\cdot \frac{2\pi i}{(\frac{n}{2})!}\cdot-\mathcal{D}_{4}dx'\nonumber\\
&+h'(0)\sum_{i=1}^{n-1}a_{i}^{i}a_{n}^{n}\cdot Vol(S^{n-2}){\rm tr}[\texttt{id}]\cdot \frac{2\pi i}{(\frac{n}{2}+1)!}\mathcal{D}_{6}dx'\nonumber\\
&-\sum_{i=1}^{n-1}g^{M}(J(e_{i}), (\nabla^{L}_{e_{i}}J)e_{n})\cdot Vol(S^{n-2}){\rm tr}[\texttt{id}]\cdot \frac{2\pi i}{(\frac{n}{2}+1)!}\cdot\mathcal{D}_{7}dx'\nonumber\\
&+\sum_{i=1}^{n-1}g^{M}(J(e_{i}), (\nabla^{L}_{e_{i}}J)e_{n})\cdot Vol(S^{n-2}){\rm tr}[\texttt{id}]\cdot \frac{2\pi i}{(\frac{n}{2}+1)!}\cdot\mathcal{D}_{8}dx'\nonumber,
\end{align}
where
\begin{align}
\mathcal{D}_{0}
&=\Big[\frac{i(n-3)\xi_n^{5}-(n^{2}-5n+8)\xi_n^{4}-i(n^{2}-3n+2)\xi_n^{3}-(3n^{2}-10n+14)\xi_n^{2}}{8\left(\xi _n+i\right)^{\frac{n}{2}+1}}\\
&+\frac{i(n^{2}-2n+1)\xi_n+n-6}{8\left(\xi _n+i\right)^{\frac{n}{2}+1}}\Big]^{(\frac{n}{2}+2)}\bigg|_{\xi_n=i};\nonumber
\end{align}
\begin{align}
\mathcal{D}_{1}
&=\Big[\frac{-i(n-3)\xi_n^{2}+i}{2\left(\xi _n+i\right)^{\frac{n}{2}}}\Big]^{(\frac{n}{2})}\bigg|_{\xi_n=i};\\
\mathcal{D}_{2}
&=\Big[\frac{in(n-3)\xi_n^{3}+2(n-3)\xi_n^{2}+i(n^{2}-7n+8)\xi_n-2n+2}{4(n-1)\left(\xi _n+i\right)^{\frac{n}{2}}}\Big]^{(\frac{n}{2}+1)}\bigg|_{\xi_n=i};\\
\mathcal{D}_{3}
&=\Big[\frac{-i(n-3)\xi_n^{4}-(n-3)\xi_n^{3}+i(n^{2}-5n+9)\xi_n^{2}+\xi_n-in+4i}{8(n-1)\left(\xi _n+i\right)^{\frac{n}{2}+1}}\Big]^{(\frac{n}{2}+1)}\bigg|_{\xi_n=i};\\
\mathcal{D}_{4}
&=\Big[-\frac{1}{8\left(\xi _n+i\right)^{\frac{n}{2}-1}}\Big]^{(\frac{n}{2})}\bigg|_{\xi_n=i};\\
\mathcal{D}_{5}
&=\Big[\frac{i(n^{2}-4n+3)\xi_n^{4}-(n^{3}-5n^{2}+5n-1)\xi_n^{3}+i(2n^{2}-6n+5)\xi_n^{2}+(n^{2}-1)\xi_n+in}{8(n-1)\left(\xi _n+i\right)^{\frac{n}{2}+1}}\Big]^{(\frac{n}{2}+1)}\bigg|_{\xi_n=i};\\
\mathcal{D}_{6}
&=\Big[\frac{(n^{2}-4n+3)\xi_n^{2}-i(n^{2}-4n+4)\xi_n-1}{8(n-1)\left(\xi _n+i\right)^{\frac{n}{2}}}\Big]^{(\frac{n}{2}+1)}\bigg|_{\xi_n=i};\\
\mathcal{D}_{7}
&=\Big[-\frac{i(n-2)\xi_n+n-2}{4(n-1)\left(\xi _n+i\right)^{\frac{n}{2}}}\Big]^{(\frac{n}{2}+1)}\bigg|_{\xi_n=i};\\
\mathcal{D}_{8}
&=\Big[\frac{(n^{2}-3n+2)\xi_n^{2}-5i(n-2)\xi_n}{4(n-1)\left(\xi _n+i\right)^{\frac{n}{2}}}\Big]^{(\frac{n}{2}+1)}\bigg|_{\xi_n=i}.
\end{align}

Through a complex series of calculations, we obtain
\begin{align}
\Phi
&=-h'(0)\cdot Vol(S^{n-2}){\rm tr}[\texttt{id}]\cdot\pi n(n^{4}-5n^{3}-16n^{2}+68n-48)2^{-4-n}\cdot\frac{(n-3)!}{(\frac{n}{2}+2)!(\frac{n}{2})!}dx'\\
&+h'(0)\sum_{\nu,i=1}^{n-1}(a_{\nu}^{i})^2\cdot Vol(S^{n-2}){\rm tr}[\texttt{id}]\cdot\pi n(n^{2}-8n+12)2^{-3-n}\cdot\frac{(n-3)!}{(\frac{n}{2}+1)!(\frac{n}{2})!}dx'\nonumber\\
&+h'(0)\sum_{i=1}^{n-1}(a_{n}^{i})^2\cdot Vol(S^{n-2}){\rm tr}[\texttt{id}]\cdot\frac{\pi(-2n^{2}+7n-2)}{n+2}2^{-n}\cdot\frac{(n-3)!}{(\frac{n}{2})!(\frac{n}{2}-2)!}dx'\nonumber\\
&-\sum_{i=1}^{n-1}g^{M}(J(e_{i}), (\nabla^{L}_{e_{i}}J)e_{n})\cdot Vol(S^{n-2}){\rm tr}[\texttt{id}]\cdot\pi(n^{2}-8n+12)2^{-1-n}\cdot\frac{(n-2)!}{(\frac{n}{2}+1)!(\frac{n}{2}-1)!}dx'\nonumber.\nonumber
\end{align}

Recall the Einstein-Hilbert action for manifolds with boundary (see \cite{Wa4} for more details)
\begin{align}
I_{G_{r}}=\frac{1}{16\pi}\int_{M}sd{\rm Vol_{M} }+2\int_{\partial M}Kd{\rm Vol_{\partial M} }:=I_{G_{r, i}}+I_{G_{r, b}},
\end{align}
where 
\begin{align}
K=\sum_{1\leq i,j\leq n-1}K_{ij}g_{\partial M}^{ij};~~~~~K_{ij}=-\Gamma_{ij}^{n},
\end{align}
and $K_{ij}$ is the second fundamental form, or extrinsic curvature. 
Take the metric of (2.1), then by \cite{Wa3}, $K_{ij}(x_{0})=-\Gamma_{ij}^{n}(x_{0})$, when $i=j<n$, otherwise it is zero. 
Then
\begin{align}
K(x_{0})&=\sum_{i,j< n-1}K_{ij}(x_{0})g_{\partial M}^{ij}(x_{0})\\
&=\sum_{i=1}^{n-1}K_{ii}(x_{0})=-\frac{n-1}{2}h'(0).\nonumber
\end{align}

A trivial verification shows that
\begin{align}
\sum_{\nu,i=1}^{n-1}(a_{\nu}^{i})^2&=\sum_{\nu,i=1}^{n}(a_{\nu}^{i})^2-\sum_{\nu=1}^{n-1}(a_{\nu}^{n})^2-\sum_{i=1}^{n-1}(a_{n}^{i})^2-(a_{n}^{n})^2\\
&=\sum_{i=1}^{n}\delta_{i}^{i}-2\sum_{i=1}^{n-1}{\langle J(e_{i}), e_{n}\rangle}^2-{\langle J(e_{n}), e_{n}\rangle}^2\nonumber\\
&=n-2\sum_{i=1}^{n}{\langle J(e_{i}), e_{n}\rangle}^2+{\langle J(e_{n}), e_{n}\rangle}^2\nonumber
\end{align}
and
\begin{align}
\sum_{i=1}^{n-1}(a_{n}^{i})^2&=\sum_{i=1}^{n}(a_{n}^{i})^2-(a_{n}^{n})^2=1-{\langle J(e_{n}), e_{n}\rangle}^2.\nonumber
\end{align}

Hence we have proved the Theorem 1.1.

\section{Declarations}
Ethics approval and consent to participate No applicable.\\

Consent for publication No applicable.\\

Availability of data and material The authors confirm that the data supporting the findings of this study are available within the article.\\

Competing interests The authors declare no conflict of interest.\\

Funding This research was funded by National Natural Science Foundation of China: No.11771070.\\

Authors' contributions All authors contributed to the study conception and design. Material preparation, data collection and analysis were performed by Siyao Liu and Yong Wang. The first draft of the manuscript was written by Siyao Liu and all authors commented on previous versions of the manuscript. All authors read and approved the final manuscript.\\

\vskip 1 true cm


\bigskip
\bigskip

\noindent {\footnotesize {\it S. Liu} \\
{School of Mathematics and Statistics, Northeast Normal University, Changchun 130024, China}\\
{Email: liusy719@nenu.edu.cn}

\noindent {\footnotesize {\it Y. Wang} \\
{School of Mathematics and Statistics, Northeast Normal University, Changchun 130024, China}\\
{Email: wangy581@nenu.edu.cn}

\end{document}